\documentclass[11pt]{article}
\usepackage{amssymb}
\usepackage{amsmath}
\usepackage{theorem}
\usepackage{latexsym}
\reversemarginpar
\theoremheaderfont{\bf} \theoremstyle{break} \theorembodyfont{\sl}
\newtheorem{theo}{Theorem}[section]
\newtheorem{defi}[theo]{Definition}
{\theorembodyfont{\rm} \newtheorem{exa}[theo]{Example}}
{\theorembodyfont{\rm} \newtheorem{rem}[theo]{Remark}}
\newtheorem{prop}[theo]{Proposition}
\newtheorem{cor}[theo]{Corollary}
\newtheorem{lemma}[theo]{Lemma}
{\theorembodyfont{\rm}\newtheorem{algorithm}[theo]{Algorithm}}
{\theorembodyfont{\rm}\newtheorem{notation}[theo]{Notation}}
\newenvironment{proof}{{\sc Proof:}}{\mbox{}\hfill$\Box$\par}
\newcommand{\eqnref}[1]{~\mbox{$(${\rm \ref{#1}}$)$}}

\newcommand{\junk}[1]{}

\newcommand{\tx}{\textstyle}
\newcommand{\N}{{\mathbb N}}
\newcommand{\F}{{\mathbb F}}

\newcommand{\cC}{{\mathcal C}}
\newcommand{\cR}{{\mathcal R}}
\newcommand{\cM}{{\mathcal M}}

\newcommand{\cJ}{{\mathcal J}}
\newcommand{\ve}[1]{{\varepsilon}^{(#1)}}
\newcommand{\LM}[1]{\mbox{\rm LM}({#1})}
\newcommand{\bi}{\begin{itemize}}
\newcommand{\ei}{\end{itemize}}
\newtheorem{obs}[theo]{Observation}

\newcommand{\bc}{\begin{center}}
\newcommand{\ec}{\end{center}}
\newcommand{\lcz}[1]{\mbox{\rm lc}_z(#1)}

\newcommand{\rank}{\mbox{\rm rank}\,}
\newcommand{\AutF}{\mbox{${\rm Aut}_{\mathbb F}$}}
\newcommand{\im}{\mbox{\rm im}\,}
\renewcommand{\mod}{\mbox{\rm mod}\,}
\newcommand{\dfree}{\mbox{$d_{\mbox{\rm\tiny free}}$}}
\newcommand{\wt}{\mbox{\rm wt}}
\newcommand{\id}{\mbox{\rm id}}

\newcommand{\ideal}[1]{\mbox{$\langle{#1}\rangle$}}

\newcommand{\Azs}{\mbox{$A[z;\sigma]$}}
\newcommand{\p}{\mbox{$\mathfrak{p}$}}
\renewcommand{\v}{\mbox{$\mathfrak{v}$}}

\newcommand{\Psigma}{\mbox{$P_{\sigma}$}}

\newcommand{\cMsigma}[1]{\mbox{$\cM^{\sigma}\!({#1})$}}
\newcommand{\trans}[1]{\mbox{$\mbox{}^{^{\scriptstyle\sf t}}\!{#1}$}}
\newcommand{\scirc}{$\sigma$-circulant}
\newcommand{\wh}[2]{\mbox{$\mbox{$\widehat{#1}$}^{\,\scriptscriptstyle{#2}}$}}
\newcommand{\rcirc}[1]{\mbox{$\mbox{${#1}$}^{^{\!\!\circ}}$}}
\newcommand{\lcirc}[1]{\mbox{$\mbox{}\hspace{.2em}\mbox{${#1}$}^{^{\hspace{-2.1em}\circ}}\hspace{2em}$}}
\newcommand{\lideal}[1]{\mbox{$^{^{\bullet\!\!}}\langle{\, #1\, }\rangle$}}
\newcommand{\rideal}[1]{\mbox{$\langle{\, #1\, }\rangle^{^{\!\!\bullet}}$}}

\newcommand{\Smallfourmat}[4]{\mbox{\tiny{$\begin{pmatrix}{#1}&{\!\!\!#2}\\{#3}&{\!\!\!#4}\end{pmatrix}$}}}

\newcounter{abc}
\newcounter{def}
\newenvironment{romanlist}{\begin{list}{(\roman{abc})\hfill}{\usecounter{abc}
     \topsep-1.4ex \labelwidth.7cm \leftmargin.7cm \labelsep0cm
     \rightmargin0cm \parsep0ex \itemsep.6ex
     \partopsep1.6ex}}{\end{list}}
\newenvironment{alphalist}{\begin{list}{(\alph{abc})\hfill}{\usecounter{abc}
     \topsep-1.4ex \labelwidth.7cm \leftmargin.7cm \labelsep0cm
     \rightmargin0cm \parsep0ex \itemsep.6ex
     \partopsep1.6ex}}{\end{list}}
\newenvironment{arabiclist}{\begin{list}{(\arabic{abc})\hfill}{\usecounter{abc}
     \topsep-1.4ex \labelwidth.7cm \leftmargin.7cm \labelsep0cm
     \rightmargin0cm \parsep0ex \itemsep.6ex
     \partopsep1.6ex}}{\end{list}}
\newenvironment{algolist}{\begin{list}{{\bf Step~\arabic{abc}:\hfill}}{\usecounter{abc}
     \topsep-1.4ex \labelwidth1.6cm \leftmargin1.6cm \labelsep0cm
     \rightmargin0cm \parsep0ex \itemsep.6ex
     \partopsep1.6ex}}{\end{list}}
\newenvironment{caselist}{\begin{list}{{\underline{Case~\arabic{def}:}\hfill}}{\usecounter{def}
     \topsep-1.4ex \labelwidth1.4cm \leftmargin1.4cm \labelsep0cm
     \rightmargin0cm \parsep0ex \itemsep.6ex
     \partopsep1.6ex}}{\end{list}}

\title{On Cyclic Convolutional Codes}
\date{\today}
\author{Heide Gluesing-Luerssen$^*$ and Wiland Schmale\footnote{
        Department of Mathematics,
        University of Oldenburg,
        26\,111 Oldenburg,
        Germany, email:
        gluesing@ mathematik.uni-oldenburg.de,
        wiland.schmale@uni-oldenburg.de}
        }
\begin{document}
\maketitle

\begin{abstract}
\noindent
We investigate the notion of cyclicity for convolutional codes as it
has been introduced in the papers~\cite{Pi76,Ro79}. Codes of this type are
described as submodules of $\F[z]^n $ with some additional
generalized cyclic structure but also as specific left ideals in a
skew polynomial ring. Extending a result of~\cite{Pi76}, we show in a
purely algebraic setting that these ideals are always principal. This
leads to the notion of a generator polynomial just like for cyclic
block codes. Similarly a control polynomial can be introduced by
considering the right annihilator ideal. An algorithmic procedure is
developed which produces unique reduced generator and control
polynomials. We also show how basic code properties and a
minimal generator matrix can be read off from these objects. A
close link between polynomial and vector description of the codes
is provided by certain generalized circulant matrices.

\end{abstract}


\section{Introduction}\label{S-intro}
\setcounter{equation}{0}
Convolutional codes (CC's) and block codes are the most widely
used types of codes in engineering practice, a fact which leads to a
continuing need for a thorough mathematical basis for the design
of useful codes. In consequence, coding theory has become one of the various
young branches of mathematics which are attractive because of the
active interplay between sophisticated engineering inventions and
high level mathematics. This is particularly true for the theory
of cyclic block codes.

The algebraic theory of CC's was initiated mainly by the articles
of Forney~\cite{Fo70} and Massey et~al.~\cite{MaSa67,MaSa68}, and, as can
be seen from the books~\cite{JoZi99,Pi88b} and the
article~\cite{McE98}, a lot of material has been accumulated
since. In the last decade Rosenthal and co-workers began a
successful project, dedicated to a better and deeper
mathematical understanding of CC's by exploiting more
systematically the existing links to control theory,
see~\cite{RSY96,RoYo99,Ro00}. Yet, up to now the mathematical
theory of CC's is not nearly as developed as that of block
codes. This gap is even larger when it comes to the notion of
cyclicity. Despite the well-known and frequently exploited efficiency of
cyclic block codes, almost nothing is known about
cyclic structures for convolutional codes and their possible impact on
applications.

In 1976 Piret showed in his fundamental paper~\cite{Pi76} how
cyclicity has to be understood for CC's and laid the basis
for a mathematical theory of cyclic CC's.
The first important discovery of Piret was that classical cyclicity ---
as common for block codes --- is trivial for CC's
(see Proposition~\ref{P-trivialCCC} in the next section).
He also showed that a more sophisticated ``graded cyclicity''
leads to interesting examples of good convolutional codes,
some of which can be found in~\cite[Sect.~IV]{Pi76} and
in~\cite{GSS02}.
His second main discovery was that irreducible graded cyclic CC's
can algebraically be described by one-sided principal ideals in a
noncommutative algebra $\Azs$.
This algebra will be introduced in the next section.
For the moment we only mention that $\Azs$ is equal to
$A[z]$ as a left $\F[z]$-module where
$A\cong \F[x]/(x^n-1)$,~$\F$ is a finite field and~$n$ is the length of
the code. Only the multiplication in the algebra $\Azs$ is quite different
from the ordinary one. It depends on an $\F$-automorphism~$\sigma$ of~$A$ and is
typically non-commutative.

The results of Piret indicate a surprising analogy to the
theory of block codes where cyclic codes are described as
principal ideals, see~\cite{MS77,Be98}, with the only difference
that the latter are in the commutative ring~$A$.

Shortly after Piret, a thorough
analysis of his results was undertaken by Roos~\cite{Ro79} in a module
theoretic framework, avoiding thereby cumbersome matrix
manipulations. At the same time Roos considerably extended
Piret's notion to what will be called $\sigma$-cyclicity later on in this paper.
But apart from this, no substantially new results could be added
and Piret's idea of a generating polynomial~\cite[Thm.3.10]{Pi76}
could not be incorporated. Furthermore, Roos' results are partly
non-constructive.

After the work of Piret and Roos no substantial effort has been made
towards a concise mathematical description of cyclic CC's  --- as far as we know.
This may partly be due to the limited mathematical readership of the
journals in question and to the circumstance that Piret's
article is quite cumbersome to read.

Yet, we think that this topic is worth being investigated in more detail.
We realized that Piret's results may serve as a good basis for a theory of
$\sigma$-cyclic CC's, which we would like to re-initiate with this paper.
Although we do not consider distance properties of cyclic CC's,
we will present the exact free distance of all codes constructed in the examples.
This way we hope to indicate that the big class of $\sigma$-cyclic CC's contains
quite some good codes and, therefore, deserves to be investigated further.

We proceed with an outline of the paper.
In Section~\ref{S-CCdef} we will trace the steps
which lead to the definition of $\sigma$-cyclicity for CC's.
We think it is worthwhile recalling also the original idea of Piret before going
into the more general setting initiated by Roos.
We will construct the (generalized) Piret algebra $\Azs$ and develop the
representation of $\sigma$-cyclic CC's as left ideals in $\Azs$.
As the Piret algebra is based on an automorphism~$\sigma$
of~$A$, we have to collect some information about the group of automorphisms of~$A$ and
about how the structure of the Piret algebra depends on~$\sigma$.
This will be done in Section~\ref{S-autos}.
In Section~\ref{S-idealgen} we give an algebraic and extended version of Piret's
main result which states that~$\sigma$-cyclic CC's are left principal ideals
in $\Azs$.
Thereafter  we investigate as to what extent a  generator of a left
ideal in $\Azs$ is unique and in Section~\ref{gencomp} we show how this unique
generator can be computed by means of a finite algorithmic procedure.
The basic algebraic tool for these sections
is a decomposition of the Piret algebra by idempotents of~$A$ and a reduction
procedure based on a monomial order of the skew polynomials.
In Section~\ref{S-circ} we introduce a new type of non-commuting $\sigma$-circulant
matrices along with a thorough investigation of their properties.
These matrices are just the proper medium for the interplay between
left ideals together with their principal generators on the one side and
CC's as submodules of $\F[z]^n$ along with their generating matrices on the other.
They also turn out to be a canonical, yet nontrivial, generalization of classical
circulants as they are common in the theory of cyclic block codes.
This becomes in particularly clear when we derive our results on generator and
control polynomials and dual codes in Section~\ref{S-gencontr}.
Indeed, we arrive at a scenario very similar to that of cyclic block codes.
The notion of a control polynomial is also included in this framework, it is obtained via
(right) annihilator ideals in the Piret algebra.
Beyond this algebraic structure, convolutional coding requires to also discuss some other
properties and invariants of the codes, as there are non-catastrophicity, minimal
generator matrices and the complexity of the given code.
All these issues can nicely be dealt with in our algebraic context.
As it turns out, the reduced principal generator polynomials for left ideals
in $\Azs$, as constructed in Section~\ref{S-idealgen} and~\ref{gencomp},
also suits well when it comes to the properties of the associate circulant matrix.
The latter leads in a canonical way to a basic minimal generator matrix of the
given code, and, consequently, the complexity can be computed in terms of the
generator polynomial.
In order to derive these results one has to combine the techniques for
circulant matrices with the algebraic methods from
Section~\ref{S-autos} --~\ref{gencomp}.
In the final Section~\ref{S-outlook} we give a short
outline of several future research topics.

Throughout this paper we make an effort to motivate and
justify the main steps of our investigations by using the classical theory of
cyclic block codes as a guideline.
We also give explicit examples in order to show how the objects in question can be computed.
This is particularly so in Section~\ref{S-gencontr} and we hope
that this way any possible impact of our results for convolutional coding can be judged more easily.

\section{What is a cyclic convolutional code?}\label{S-CCdef}
\setcounter{equation}{0}

In this section we will shortly recall the basic definitions and
properties of convolutional codes and cyclic block codes and then
develop --- along the lines of the articles~\cite{Pi76,Ro79} ---
the notion of a cyclic convolutional code.

Throughout this paper, $\F$ denotes a fixed finite field and~$n$ a positive
integer such that
\begin{equation}\label{e-coprime}
  \text{the characteristic of~$\F$ does not divide $n$.}
\end{equation}
The number~$n$ is going  to be the length of the code and\eqnref{e-coprime}
is the familiar assumption from the theory of cyclic block codes guaranteeing that
the polynomial $x^n-1$ factors into different prime polynomials over~$\F$.

As is well-known, a block code is simply a subspace of the vector space~$\F^n$.
Analogously, convolutional codes are direct summands of $\F[z]^n$.
Of course, only additional properties single out the codes which are
relevant for applications.
Before presenting the according notions, we first collect some basic facts about
submodules and direct summands of $\F[z]^n$.

As usual in coding theory, all vectors are regarded as row vectors, thus
\[
   \F[z]^n:=\{(v_1,\ldots,v_n)\mid v_i\in\F[z]\text{ for }i=1,\ldots,n\}.
\]
Consequently, images and kernels of matrices will always denote {\em left\/} images
and {\em left\/} kernels.
The following facts will be used freely.

\begin{prop}\label{P-basicLA}
Let $V$ be a submodule of $\F[z]^n$.
\begin{alphalist}
\item $V$ has a finite basis and all bases of~$V$ have the same length, called the
      rank of~$V$.
\item If $v_1,\ldots,v_r\in\F[z]^n$ form a generating set of~$V$, then
      $V=\im M:=\{uM\mid u\in\F[z]^r\}$ where
      \begin{equation}\label{e-genmat}
          M:=\begin{bmatrix}v_1\\ \vdots\\v_r\end{bmatrix}\in\F[z]^{r\times n}.
      \end{equation}
      We call~$M$ a generating matrix of~$V$.
\item Let $P\in\F[z]^{r\times r}$ and~$M$ as in~(b). Then
      $V=\im(PM)\Longleftrightarrow P\text{ is invertible over }\F[z]$.
\end{alphalist}
\end{prop}

The following properties about direct summands are easily obtained from
linear algebra over the Euclidean domain~$\F[z]$.

\begin{prop}\label{P-directsummand}
Let $V\subseteq\F[z]^n$ be a submodule and $v_1,\ldots,v_r\in\F[z]^n$ a generating
set of~$V$. Put $M\in\F[z]^{r\times n}$ as in\eqnref{e-genmat}.
Then the following are equivalent.
\begin{arabiclist}
\item $V$ is a direct summand of $\F[z]^n$.
\item Any basis of~$V$ can be completed to a basis of $\F[z]^n$.
\item The Smith-form of~$M$ is given by
      $\Smallfourmat{I_k}{0}{0}{0}$, where~$k$ is the rank of~$V$.
\item If $v_1,\ldots,v_r$ is a basis of~$V$ (equivalently, if~$r$ is the rank
      of~$V$), then~$M$ is right invertible over $\F[z]$.
\item For all $v\in V$ and all $\lambda\in\F[z]\backslash\{0\}$ one has
      \begin{equation}\label{e-lambdav}
           \lambda v\in V\Longrightarrow v\in V.
      \end{equation}
\item There exists some matrix $N\in\F[z]^{n\times l}$ such that
      $V=\ker N:=\{v\in\F[z]^n\mid vN=0\}$.
\item For all submodules $W\in\F[z]^n$ having the same rank as~$V$ one has
      \[
          V\subseteq W\Longrightarrow V=W.
      \]
\end{arabiclist}
A matrix~$M$ with property~(3) will be called basic.
\end{prop}

For the definition of Smith-forms see e.~g. \cite[p.~141]{Ga77}
or~\cite[Sec.~3.7]{Ja85}.

A convolutional code is simply defined to be a direct summand of $\F[z]^n$.
But of course only various additional notions lead to useful coding theoretical
concepts.

\begin{defi}\label{D-CC}
\begin{arabiclist}
\item A convolutional code (CC) of length~$n$ and dimension~$k$ is a direct summand~$\cC$ of
      $\F[z]^n$of rank~$k$.
\item Any generating matrix $G\in\F[z]^{k\times n}$ having rank~$k$ of~$\cC$ is called a
      generator matrix or encoder of~$\cC$. Hence
      \[
        \cC=\im G=\{uG\mid u\in\F[z]^k\}
      \]
      and the vector $uG\in\F[z]^n$ is said to be the codeword associated with the
      message word $u\in\F[z]^k$.
\item A matrix~$H\in\F[z]^{n\times(n-k)}$ satisfying
      $\cC=\ker H=\{v\in\F^n:\, vH=0\} $ is said to be a
      control matrix of the code~$\cC$.
\item The maximal degree of the $k$-minors of an encoder~$G$ is called
      the complexity  of the code.
      A code of complexity zero is said to be a block code.
\end{arabiclist}
\end{defi}

Notice that each code has a generator and a control matrix. The control
matrix always has rank~$n-k$.
We would like to point out the difference between a generator matrix and a
generating matrix in the sense of Proposition~\ref{P-basicLA}:
the latter one need
not have full rank and therefore is not suitable as an encoder.
However, we will need this notion, since certain square (singular) generating
matrices naturally arise in our investigations of cyclic convolutional codes.
Of course, one can always constructively obtain a (full rank) generator matrix out
of these matrices by computing for instance the Hermite normal form.
But in our specific context a better way to a generator matrix will be shown in
Section~\ref{S-gencontr}.

\begin{rem}\label{R-CC}
\begin{arabiclist}
\item In Definition \ref{D-CC} we adopt the viewpoint that codewords and message words are finite
      sequences rather than infinite ones, the latter being slightly more common in
      convolutional coding theory; see~\cite{RSY96} for a discussion of this subtle difference.
      The codewords and messages are therefore represented by polynomials rather than by Laurent series
      from $\F(\!(z)\!)$.
      But in any case, even if Laurent series are admitted,
      the encoders are always polynomial matrices exactly as in Definition~\ref{D-CC}, see,
      e.~g., \cite{Fo70,McE98}.
      Moreover, there is a one-one correspondence between CC's in the
      sense of Definition~\ref{D-CC} and CC's as subspaces of
      $\F(\!(z)\!)^n$ with a polynomial generator matrix.
\item It is well-known that the complexity does not depend on the choice of the
      encoder.
      Furthermore, from the theory of minimal bases (see~\cite{Fo75}) it follows
      that a code has complexity zero if and only if it has a constant encoder.
      Thus, in this case the code behaves just like a block code.
\end{arabiclist}
\end{rem}

Some of our investigations will hold under weaker assumptions,
which are closely related to the following concepts of coding theory.

\begin{defi}\label{D-attrib}
Let $V$ be a submodule of $\F[z]^n$. Then
\begin{alphalist}
\item $V$ is called non-catastrophic if\eqnref{e-lambdav}
      is satisfied for all $v\in V$ and all $\lambda\in\F[z]\backslash z\F[z]$.
\item $V$ is called delay-free if\eqnref{e-lambdav}
      is satisfied for all $v\in V$ and $\lambda=z$.
\end{alphalist}
\end{defi}

A non-catastrophic and delay-free submodule is, by definition, a
convolutional code, see Proposition~\ref{P-directsummand}(5).

A first indication for the quality of a code is given by its free
distance defined as follows, for details see also \cite[Ch.~3]{JoZi99}.

\begin{defi}\label{D-dist}
\begin{alphalist}
\item The weight of the word\, $v=\sum_{\nu=0}^N v_{\nu}z^{\nu}\in\F[z]^n$ is defined
      as $\wt(v):=\sum_{\nu=0}^N\wt(v_{\nu})$ $\in\N_0$,
      where $\wt(v_{\nu})$ denotes the Hamming weight of the constant
      vector~$v_{\nu}\in\F^n$.
\item The free distance of a code $\cC\subseteq\F[z]^n$ is defined as
      $\dfree(\cC):=\min\{\wt(v)\mid v\in\cC\backslash\{0\}\}$.
\end{alphalist}
\end{defi}

In our examples we will often state explicitly the free distance of the code under investigation.
Although we do not investigate the free distance of a cyclic convolutional code in
this paper, we think it worthwhile computing the distance in
order to have a more complete picture of the codes in question.
Most of these computations have been done with the help of MAPLE.

Now we turn to the notion of cyclicity. A block code
$\cC\subseteq\F^n$ is said to be cyclic if it is invariant under the
cyclic shift, that is, if
\begin{equation}\label{e-shift}
  (v_0,\ldots,v_{n-1})\in\cC\Longrightarrow
  (v_{n-1},v_0,\ldots,v_{n-2})\in\cC
\end{equation}
or, equivalently, if
\begin{equation}\label{e-CS}
   \cC S\subseteq\cC
\end{equation}
where
\begin{equation}\label{e-shift-2}
   S=\begin{bmatrix}0&1&\cdots  & 0 \\
     \vdots &\vdots &\ddots&\vdots  \\
     0 & 0 &\cdots  & 1\\ 1& 0 &\cdots &0
     \end{bmatrix}\in\F^{n\times n}.
\end{equation}

An important tool in the theory of cyclic block codes is the
so-called polynomial representation. It is based on the
$\F$-isomorphism
\[
  \p: \F^n\longrightarrow A,\quad v=(v_0,\ldots,v_{n-1})\longmapsto
      \p (v)=\sum_{i=0}^{n-1}v_ix^i,
\]
where $A:=\F[x]/_{\ideal{x^n-1}}$ is displayed in the
canonical way
\[
A=\{f\in\F[x]\mid \text{deg}f<n\}\text{ with multiplication modulo }
x^n-1.
\]
The inverse of $\p$  will be denoted by $\v$. The map~$\p$ translates
the cyclic shift into multiplication by~$x$. As a consequence, a cyclic block code
$\cC$  can now be represented as an ideal $\p(\cC)$ in $A$ and vice versa; in other
words,
\begin{equation}\label{E-blockcyc}
    \text{a block code $\cC$ is cyclic if and only if}\
    \big[\,a\in\p(\cC) \Longrightarrow\ xa\in\p(\cC)\,\big].
\end{equation}
For later use we immediately extend $\p$ to all of $\F[z]^n$ via
\begin{equation}\label{E-p}
      \p:\F[z]^n\, \longrightarrow\, A[z],\quad
        \sum_{\nu\geq0}z^{\nu}v_{\nu}
        \longmapsto\sum_{\nu\geq0}z^{\nu}\p (v_{\nu}).
\end{equation}
The map $\p$ is an isomorphism of left $\F[z]$-modules with inverse
$\v:=\p^{-1}$.

It would be quite natural to define cyclicity of convolutional
codes just like for block codes, that is, by requiring invariance
as in\eqnref{e-shift}. But already in \cite[Thm.~3.12]{Pi76} and
\cite[Thm.6]{Ro79} the following important observation has been made.

\begin{prop}\label{P-trivialCCC}
Let $\cC\subseteq\F[z]^n$ be a code satisfying\eqnref{e-shift}.
Then~$\cC$ is a block code.
\end{prop}

This result will appear as a special case in Proposition \ref{P-blockcode}.
However, we include an independent and elementary linear algebraic proof at the end
of this section.

The negative result of Proposition \ref{P-trivialCCC} has led
Piret~\cite{Pi76} to a more general and complex notion of cyclicity for convolutional codes.
Instead of shift-invariance of~$\cC$ under the shift-matrix~$S$ from\eqnref{e-shift-2},
which would require $\sum_{{\nu}=0}^d z^{\nu} v_{\nu} S\in\cC$, whenever
$\sum_{{\nu}=0}^d z^{\nu} v_{\nu} \in\cC$,
Piret introduced a kind of graded quasi-cyclicity. Precisely, he called a convolutional
code~$\cC$ cyclic, if there exists some~$m$, which is coprime to the length~$n$ of the code,
such that
\begin{equation}\label{Piret-def1}
   \sum z^{\nu} v_{\nu} \in\cC \, \Longrightarrow\,
   \sum z^{\nu} v_{\nu} S^{(m^{\nu})} \in\cC.
\end{equation}
In polynomial language, i.~e. in the polynomial ring~$A[z]$, this translates into
\begin{equation}\label{Piret-def2}
   \sum z^{\nu} \p (v_{\nu}) \in\p (\cC)\,
   \Longrightarrow\,
   \sum z^{\nu} x^{(m^{\nu})} \p(v_{\nu}) \in \p(\cC).
\end{equation}
The coprimeness of the integers~$m$ and~$n$ guarantees not only that the
minimal polynomial of $S^m$ is the same as that of~$S$, that is $x^n-1$,
but also that the map $x \longmapsto x^m$ induces an $\F$-automorphism of $A$.
This allows to introduce an $\F$-algebra structure on the left
$\F[z]$-module $A[z]$ which naturally extends the algebra A. The
details of the construction will be explained below.

Piret's notion of cyclicity has been generalized by Roos~\cite{Ro79} in a natural way to
arbitrary $\F$-automorphisms~$\sigma$ of~$A$.
We propose the name
{\it $\sigma$-cyclicity}, since later on different automorphisms
will have to be considered simultaneously.
In the following definition we introduce this notion for
arbitrary submodules of $\F[z]^n$.

\begin{defi}\label{D-Roos}
Let $\AutF(A)$ denote the group of all~$\F$-algebra automorphisms on
$A$ and let $\sigma\in\AutF(A)$.
A submodule $\cC$ of $\F[z]^n$ is called $\sigma$-cyclic if
\begin{equation}\label{E-Roos}
   g=\sum_{\nu=0}^dz^{\nu}g_{\nu}\in\p(\cC)\,\Longrightarrow\, x*_{\sigma}g:=
   \sum_{\nu=0}^dz^{\nu}\sigma^{\nu}(x)g_{\nu}
   \in\p(\cC).
\end{equation}
Consequently, a $\sigma$-cyclic convolutional code ($\sigma$-CCC) is a $\sigma$-cyclic
direct summand of $\F[z]^n$.
\end{defi}

In~\cite{Ro79}, Equation\eqnref{E-Roos} was extended to a left $\F[z]$-module
structure on~$A[z]$, which then was used to investigate in great detail the
structure of $\sigma$-CCCs.
Unfortunately, generator polynomials as constructed by Piret could not be
incorporated in this setting.
It seems to be more helpful to use~$*_{\sigma}$ for a non-commutative
ring structure on $A[z]$ as follows.

\begin{defi}\label{D-prod}
Let $\sigma\in\AutF(A)$. We define the product of
$g=\sum_{\nu\ge 0}z^{\nu}g_{\nu},\;h=\sum_{\mu\ge 0}z^{\mu}h_{\mu}\in A[z]$ by
\[
  \big(\sum_{\nu\geq0} z^{\nu}g_{\nu}\big)*_{\sigma}
  \big(\sum_{\mu\geq0} z^{\mu}h_{\mu}\big)
  :=\sum_{\lambda\geq0}z^{\lambda}\sum_{\nu+\mu=
  \lambda}\sigma^{\mu}(g_{\nu})h_{\mu}.
\]
$A[z]$ equipped with the multiplication~$*_{\sigma}$ will be denoted by
$\Azs$ and often be abbreviated by $\cR$.
We call $\Azs$ a Piret algebra (with parameters $q=|\F|,\,n,\,\sigma$).
\end{defi}

Observe that multiplication in $\Azs$ is simply an extension of the (commutative)
multiplication in~$A$ together with the rule
\begin{equation}\label{e-az}
   a*_{\sigma}z=z*_{\sigma}\sigma(a)\text{ for all }a\in A.
\end{equation}
In particular we have
\[
   \lambda*_{\sigma}z=z*_{\sigma}\lambda\text{ for all }\lambda\in\F
\]
and therefore we obtain the usual product whenever the left factor is in $\F[z]$;
precisely,
\[
  \big(\sum_{\nu\geq0} z^{\nu}g_{\nu}\big)*_{\sigma}
  \big(\sum_{\mu\geq0} z^{\mu}h_{\mu}\big)
  =\sum_{\lambda\geq0}z^{\lambda}\sum_{\nu+\mu=
  \lambda}g_{\nu}h_{\mu}\text{ for all }
  \sum_{\nu\geq0} z^{\nu}g_{\nu}\in\F[z]\text{ and }
  \sum_{\mu\geq0} z^{\mu}h_{\mu}\in A[z].
\]
Notice that we put the $z$-coefficients always to the right of~$z$.
This is, of course, a matter of choice, but the explicit form of the
non-commutative product $*_{\sigma}$ highly depends on it. Using
the multiplication rule\eqnref{e-az} the coefficients can always be
shifted to the left of~$z$ if a suitable power of $\sigma^{-1}$ is
applied.
One should always bear in mind that a monomial $z^{\nu}a,\,a\in A$,
can also be read as $z^{\nu}*_{\sigma}a$.

The notation $\Azs$ is common in the theory of skew polynomial rings over integral domains,
see for instance~\cite[p.~438]{Co77}.
In our setting~$\Azs$ typically has many zero divisors and --- as
we will see later --- many nonconstant units.

The discussion above leads directly to a nice skew polynomial
representation for $\sigma $-cyclic submodules.

\begin{obs}\label{Piretalgebra}
\begin{alphalist}
\item $\Azs$ is an $\F$-algebra which, at the same time, carries a canonical left
      $\F[z]$-module structure.
      The algebra $\Azs$ is non-commutative whenever the automorphism~$\sigma$ is not the identity
      on~$A$.
\item A submodule~$\cC$ of~$\F[z]^n$ is $\sigma$-cyclic if and only if its polynomial
      version $\p(\cC)$ is a left ideal in $\Azs$.
\item It is worthwhile noting that $\Azs$ also carries a canonical right $\F[z]$-structure.
      This is even more obvious than the left $\F[z]$-module structure since in the
      multiplication $*_{\sigma}$ the automorphism~$\sigma$ acts only on the left
      factor.
\end{alphalist}
\end{obs}

\begin{exa}\label{E-CCC1}
Let $\F=\F_4=\{0,1,\alpha,\alpha^2\}$ and $n=3$.
\begin{arabiclist}
\item  We choose the automorphism~$\sigma$ given by $\sigma(x)=\alpha^2 x$
       (it will be explained in Example~\ref{E-Azs1} below that this indeed induces an
       automorphism on~$A$).
       We wish to find the smallest $\sigma$-CCC~$\cC$ containing the codeword
       \[
         v:=(1+z+z^2,\;\alpha+z+\alpha^2 z^2,\;\alpha^2+z+\alpha z^2)\in\F[z]^3.
       \]
       First of all, $\p(\cC)$ has to contain the left ideal in $\Azs$ generated by
       the polynomial
       \[
          g:=\p(v)=1+\alpha x+\alpha^2 x^2+z(1+x+x^2)+z^2(1+\alpha^2 x+\alpha x^2).
       \]
       One calculates
       \[
         x*_{\sigma}g=\alpha^2+x+\alpha x^2+z\alpha^2(1+x+x^2)+z^2(\alpha^2+\alpha x+x^2)
           =\alpha^2 g
       \]
       and thus $x^2*_{\sigma}g=\alpha g$. Furthermore, one easily checks that the matrix
       \[
         G:=\big[1+z+z^2,\;\alpha+z+\alpha^2 z^2,\;\alpha^2+z+\alpha z^2\big]
       \]
       is basic and therefore $\cC=\im G$ is the smallest $\sigma$-CCC containing
       the word~$v$ above.
       This code happens to be quite a good one, since one
       can show that $\dfree(\cC)=9$, which is the maximum value for the free
       distance of any one-dimensional code of length~$3$ and complexity~$2$,
       see~\cite[Thm.~2.2]{RoSm99}.
       Hence~$\cC$ is an MDS-code in the sense of~\cite[Def.~2.5]{RoSm99}.
\item Let us also consider the situation in~(1) with the automorphism
      $\sigma=\mbox{\rm id}$.
      In this case multiplication by~$x$ simply corresponds to the usual cyclic shift and
      therefore the smallest $\sigma$-CCC $\cC'$ containing~$v$ has to satisfy
      \[
        \cC'\supseteq\im G',\text{ where }
        G':=\begin{bmatrix}
               1+z+z^2&\alpha+z+\alpha^2 z^2&\alpha^2+z+\alpha z^2\\
               \alpha^2+z+\alpha z^2&1+z+z^2&\alpha+z+\alpha^2 z^2\\
               \alpha+z+\alpha^2 z^2&\alpha^2+z+\alpha z^2&1+z+z^2\end{bmatrix}.
      \]
      Since $\det G'\not=0$, the code $\cC'$ is $3$-dimensional and, by
      Proposition~\ref{P-directsummand}(7), it follows
      \[
        \cC'=\im I_3=\F[z]^3.
      \]
      Hence $\cC'$ is a (trivial) block code and we encounter an example of the
      result in Proposition~\ref{P-trivialCCC}.
\item In the paper~\cite{Pi88} Piret gave a class of unit memory convolutional
      codes based on Reed-Solomon block codes.
      One can show that these codes are all $\sigma$-cyclic with respect to the
      automorphism given by $\sigma(x)=x^{n-1}$.
\item Finally, we would like to mention the class of convolutional
      codes constructed in the paper~\cite{SGR01}.
      Just like the codes in~\cite{Pi88} they are based on cyclic block codes and, therefore, have a
      generator matrix with a type of row-wise cyclic shift structure.
      Yet, they are in general {\em not\/} $\sigma$-cyclic with respect to any
      automorphism~$\sigma$.
      \mbox{}\hfill$\Box$
\end{arabiclist}
\end{exa}

As has been explained above, Definition \ref{D-prod} and Observation \ref{Piretalgebra}
basically go back to~\cite{Pi76}, with the only difference that in~\cite{Pi76} only monomial
automorphisms are considered, i.~e. automorphisms~$\sigma$, where $\sigma(x)=x^m$ for
some $m\in\N$.
It is easy to see that the set $\{m\mid 1\leq m\leq n-1,\,\gcd(m,n)=1\}$ leads to all
monomial automorphisms.
Note also that for every~$n$, the choice $m=n-1$ produces the automorphism given by
$\sigma (x)=x^{-1}$.

\begin{rem}\label{R-CCC}
Definition \ref{D-Roos} extends cyclicity of block codes in the
sense of\eqnref{E-blockcyc}.
One can also express $\sigma$-cyclicity
solely in terms of vector polynomials, i.~e., without resorting to the
identifications~$\p$ and~$\v$.
This yields a generalization of cyclic block codes in the sense of\eqnref{e-shift}.
Since this is more easily understood after some appropriate objects have been
defined, we will postpone this description to Observation~\ref{P-cMsigma2}.
At this point one should simply note that for
$\sigma=\mbox{\rm id}$ one has
\[
   \v\big(x*_{\sigma}\p(v)\big)=vS\quad \text{ for } \ v\in\F[z]^n
\]
(the usual cyclic shift) and in this case the map $v\mapsto\v\big(x*_{\sigma}\p(v)\big)$ on
$\F[z]^n$ is $\F[z]$-linear.
For $\sigma\not=\mbox{\rm id}$ this is no longer true, due to
non-commutativity of $\Azs$.
\end{rem}

\begin{exa}\label{E-Azs1}
The above raises the question as to how the group $\AutF(A)$ looks
like. A very simple, but tedious way of finding all automorphisms is
as  follows. First of all, notice that any $\F$-algebra
automorphism~$\sigma$ is fully determined by the value of $\sigma(x)$
in~$A$. Secondly, since~$x^n=1$ and
$1,x,\ldots,x^{n-1}$ are linearly independent over~$\F$, the same has to
be true for $a:=\sigma(x)\in A$. Furthermore, it is easy to see that
each element $a\in A$ such that $1,a,\ldots,a^{n-1}$ are linearly
independent and $a^n=1$ uniquely determines an automorphism
$\sigma\in\AutF(A)$ via
$\sigma(x)=a$. Of course, $a=x$ corresponds to $\sigma=\mbox{\rm id}$.
Applying this for instance to the case
$\F=\F_4=\{0,1,\alpha,\alpha^2\}$ and $n=3$ leads to six
automorphisms given by
\[
  a\in\{x, x^2, \alpha x, \alpha^2x, \alpha x^2, \alpha^2 x^2\}.
\]
In the next section a more sophisticated and detailed investigation
of the group $\AutF(A)$ will be presented.
\\
For an example of the non-commutativity of $\Azs$ take
e.~g. the isomorphism~$\sigma$ given by $\sigma(x)=\alpha x$.
Then $x^2*_{\sigma}z=z*_{\sigma}\sigma(x^2)=z*_{\sigma}\alpha^2x^2$.
\mbox{}\hfill$\Box$
\end{exa}

In the rest of this paper we will omit the symbol $*_{\sigma}$ in the skew
multiplication of Definition~\ref{D-prod}.
Precisely,
\begin{equation}\label{e-skewprod}
    gh:=g*_{\sigma}h\text{ for all }g,\,h\in\Azs.
\end{equation}
This won't cause any confusion since the Piret-algebra under investigation will always be
clear from the context.

Since we will often switch between $\sigma$-cyclic submodules of $\F[z]^n$ and their
counterpart as left ideals in the Piret-algebra, the following will be very
convenient. Notice that we make use of the notation in\eqnref{e-skewprod}.

\begin{obs}\label{O-attrib}
Let $\sigma\in\AutF(A)$. A left ideal~$\cJ$ of $\Azs$
is called non-catastrophic (resp.\ delay-free) if $\v(\cJ)$ is a non-catastrophic
(resp.\ delay-free) submodule of $\F[z]^n$.
Since~$\p$ and~$\v$ are $\F[z]$-linear mappings, this is equivalent to
\\[.5ex]
$
\mbox{}\!\!\begin{array}{l}
  \cJ\text{ non-catastrophic }\!\!\Longleftrightarrow\!\! \big[\,\forall\; g\in\Azs\ \forall\;
  \lambda\in\F[z]\setminus z\F[z]:\; \lambda g\in\cJ\Longrightarrow
  g\in\cJ\,\big],\\[.5ex]
  \cJ\text{ delay-free }\!\!\Longleftrightarrow\!\!\big[\, \forall\; g\in\Azs\ \forall\; k\ge 1
    :\; z^k g\in\cJ\Longrightarrow g\in\cJ\,\big],\\[.5ex]
  \cJ\text{ is a direct summand of }\Azs\!\!\Longleftrightarrow\!\!\big[\, \forall\, g\in\Azs\;\forall\;
  \lambda\in\F[z]\backslash\{0\}\!:\, \lambda g\in\cJ\Longrightarrow g\in\cJ\,\big],
\end{array}
$
\\[.5ex]
where a direct summand is understood in the context of left $\F[z]$-modules. Recall
that~$\cJ$ is a direct summand if and only if $\v(\cJ)$ is a convolutional code.
\\
We will also need the corresponding notions for right ideals, in which case, of
course, $\lambda g$ and $z^kg$ have to be replaced by $g\lambda$ and
$gz^k$, respectively. In this case, one has to recall from Observation~\ref{Piretalgebra}
that $\Azs$ is also a right $\F[z]$-module.
\end{obs}

We conclude this section with a direct proof
of Proposition \ref{P-trivialCCC}.

\begin{proof}
By assumption $\cC S\subseteq \cC$, where~$S$ is as in\eqnref{e-shift-2}.
The minimal polynomial of~$S$ is given by $x^n-1$.
Let $x^n-1=\pi_1\cdots \pi_r $ be the factorization into prime polynomials, which
are, due to\eqnref{e-coprime}, pairwise different.
Then we obtain the decomposition
\[
   \F[z]^n=\ker \pi_1(S) \oplus \cdots \oplus \ker \pi_r(S)
\]
of $\F[z]^n$ into $\F[z]$-submodules which are minimal $S$-invariant
direct summands.
Since~$\cC$ itself is a direct summand, too, we similarly obtain
\[
   \cC= \bigoplus_{i\in T}\ker \pi_i(S),\ {\rm where}\ T=\{i \mid \ker
   \pi_i(S)\cap\cC\ne \{0\}\}.
\]
Since $\F^n S=\F^n $, the $\F[z]$-submodules
$ \ker \pi_i(S)$ are generated by $\ker \pi_i(S)\cap \F^n$ and
this leads directly to a constant generating matrix and
thus to a constant encoder for~$\cC$.
By Definition~\ref{D-CC}(4) the complexity is zero, i.~e.~$\cC$ is a block code .
\end{proof}

\section{Basic information on $\F $-automorphisms            
of \protect{$\Azs$}}\label{S-autos} \setcounter{equation}{0} 

As is clear from the last section, in order to get access to all
$\sigma$-cyclic convolutional codes, it is necessary to have precise
information on the group $\AutF(A)$ and its action on the
components of $A=\F[x]/\ideal{x^n-1}$ when represented as a cartesian
product of fields. We now give this information as far as absolutely
necessary  and for reasons  of space only partially with proofs.

Under the assumption\eqnref{e-coprime} we know that the normalized
factors $\pi_i\in\F[x]$ of the prime factor decomposition
\begin{equation}\label{E-xn-1}
   x^n-1=\pi_1\cdot\ldots\cdot\pi_r
\end{equation}
are pairwise different.
We order the prime polynomials such that
$$ \deg\pi_1=\dots =\deg\pi_{r_1} < \cdots \cdots
 < \deg\pi_{r_1+\dots + r_{s-1}+1}=\dots =\deg\pi_{r_1+\dots + r_s}$$
where $r_1+ \dots + r_s=r$.

The most natural and constructive way to represent and decompose the
$\F$-algebra $A$ is as follows. Let
\begin{equation}\label{E-A}
A:=\{f\in \F[x]\mid \deg f < n \} \text{  with multiplication modulo
} x^n-1
\end{equation} and for $1\leq k\leq r$ let
\begin{equation}\label{E-Ki}
K_k:=\{f\in \F[x]\mid \deg f < \deg \pi_k \} \text{  with
multiplication modulo } \pi_k.
\end{equation}
$K_k\cong\F[x]/\ideal{\pi_k}$ is a finite Galois extension of $\F$ of dimension
$[K_k:\F]=\deg\pi_k$.
Denote by $\varrho_k(a)\in K_k$ the remainder of $a\in\F[x]$ when
dividing by $\pi_k $. By means of the Chinese remainder theorem the
map
\begin{equation}\label{E-rho}
 \varrho : A  \longrightarrow K_1 \times \dots \times K_r, \quad
 a\longmapsto[\varrho_1 (a),\cdots , \varrho_r (a)]
\end{equation}
is an isomorphism of rings, where the cartesian product is endowed
with component-wise addition and multiplication. The isomorphism
$\varrho $ can be computed easily and it induces an isomorphism of
the respective automorphism groups. Therefore, in this section we
assume from now on that
\begin{equation}\label{E-Adecomp}
A=K_1 \times \ldots \times K_r.
\end{equation}
The basic properties of the ring~$A$ which we will use in the
following reflect the fact that~$A$ is a semi-simple ring. The
canonical $\F$-basis vectors
\begin{equation}\label{E-canIdem}
   \ve{k}=[\delta_{k,j}]_{1\le j \le r}=[0,\dots,0,1,0,\dots,0],
   \ \text{where the $1$ is at the $i$-th position},
\end{equation}
are at the same time the uniquely determined primitive and pairwise
orthogonal idempotents of $A$. Recall that an idempotent is called
{\em primitive\/} if it cannot be written as a nontrivial sum of
orthogonal idempotents. We call
\begin{equation}\label{E-Kk}
K^{(k)}:= \ve{k}A=0\times \cdots \times K_k \times \cdots \times 0
\end{equation}
the $k$-th component of~$A$. Each component of~$A$ is a field, since
of course $K^{(k)}\cong K_k$.
In particular one has for all $a,\,b\in A$ the rule
\begin{equation}\label{E-idemzero}
   a\ve{k}b=0\Longrightarrow a\ve{k}=0 \text{ or } \ve{k}b=0.
\end{equation}
Any ideal of $A$ is readily seen to be of the type
$$ \sum_{k=1}^r U_k \text{ where } U_k\in \{
\{0\},K^{(k)}\} \text{ for } 1\le k \le r.$$ Two components $K^{(k)}$
and $K^{(l)}$ are isomorphic if and only if
$\deg \pi_k=\deg \pi_l$.
Therefore up to a further, usually non-unique, automorphism we can
even assume from now on that
\begin{equation}\label{E-Ldecomp}
 A= \underbrace{L_1\times \cdots  \times L_1}_{r_1}\times \cdots\cdots \times
 \underbrace{L_s\times \cdots  \times L_s}_{r_s}
 =L_1^{r_1}\times \cdots \times L_s^{r_s}
 \end{equation}
where the field $L_j$ is isomorphic to
$K_{\sum_{k=1}^{j-1}r_k+1}\cong\cdots\cong
 K_{\sum_{k=1}^j r_k}$ and, as a consequence,
$L_1,\ldots,L_s$ are pairwise non-isomorphic.
In particular $\sum_{k=1}^s r_k=r$.

Let us now consider the $\F$-automorphisms of~$A$. One first observes
that for an automorphism $\sigma\in\AutF(A)$ necessarily
$\sigma(K^{(k)})=K^{(l)}$ for some $l$. Thus $\sigma$ acts as a permutation on
the set $\mathcal K =\{K^{(1)},\dots ,K^{(r)}\}$ of components of $A$
and $\mathcal K$ is the disjoint union of cycles determined by
$\sigma$. All fields in one cycle  must have the same degree over
$\F$ and therefore are isomorphic. Therefore~$\sigma$ can only
permute those components of~$A$ which correspond to one of the
fields~$L_j$ for a fixed~$j$ in the decomposition\eqnref{E-Ldecomp}. On the
other hand, any such type of permutation together with automorphisms
of the components induces an $\F$-automorphism of~$A$ and it can be
shown that there are no further automorphisms. This is the main
information of the following fundamental theorem.

\begin{theo}\label{autos}
For $1\le j\le s$ let $G_j:=\AutF (L_j)$. Let furthermore
$S_{r_1,\dots ,r_s}$ be the subgroup of the group $S_r$  of
permutations of $\{1,\dots ,r\}$, which leaves all sets
\[
\{1,\dots ,r_1\}\ ,\ldots \ ,\  \{\,{\tx\sum_{k=1}^{s-1}r_k+1},\dots
,
 \underbrace{{\tx\sum_{k=1}^s r_k}}_{=r} \}
\]
invariant. Then
\begin{equation}\label{wreath}
   \AutF (A)\cong (G_1^{r_1}\times \cdots \times G_s^{r_s})
                  \circ(S_{r_1,\dots ,r_s}),
\end{equation}
where  $\circ$ is defined as
$$
 \big((\gamma_1,\dots ,\gamma_r)\circ\beta\big)[a_1,\dots ,a_r]=
 \big[\gamma_{\beta(1)}(a_{\beta(1)}),\dots,\,
 \gamma_{\beta(r)}(a_{\beta(r)})\big]$$
for all $[a_1,\dots,a_r]\in L_1^{r_1}\times \cdots \times L_s^{r_s}$.
\end{theo}

Note that the group on the right hand side of\eqnref{wreath} is the
automorphism group of~$A$ in the representation\eqnref{E-Ldecomp} and
only upon incorporating a fixed isomorphism leading from\eqnref{E-A}
to\eqnref{E-Ldecomp} one obtains the isomorphisms for~$A$ in the
description of\eqnref{E-A}. We will describe this translation in
detail via an example below. The representation in\eqnref{wreath} is
an instance of the wreath product. In~\cite{Ve85} one can find in a
more general situation a result (without proof) from which
Theorem~\ref{autos} could be deduced. For our purposes a direct proof
of the Theorem is preferable and not very difficult in the concrete
context as developed before Theorem~\ref{autos}. However, we skip the
proof for the sake of brevity. As an immediate consequence we obtain
a formula for the number of automorphisms on~$A$.

\begin{cor}\label{autnumb}
Let the data be as in\eqnref{E-xn-1} and\eqnref{E-Ldecomp}. Then
$|\AutF(A)|=(\deg\pi_1)^{r_1}\cdots(\deg\pi_s)^{r_s}
  r_1!\, \cdots \, r_s!$.
\end{cor}

The advantage of Theorem~\ref{autos} is that it provides us with a
very systematic and well-organized list of the automorphisms on~$A$
in the representation\eqnref{E-Ldecomp}. However, for the
investigations of cyclic codes in Section~\ref{S-circ} and
thereafter, we will need the $\F$-automorphisms for the ring~$A$ as
given in\eqnref{E-A}, i.~e. for $\sigma\in\AutF(A)$ we will need to
know the value $\sigma(x)\in A$, which completely
determines~$\sigma$. In order to find this representation of $\sigma$
one has to incorporate an isomorphism leading from\eqnref{E-A}
to\eqnref{E-Ldecomp}. This is illustrated in
Example~\ref{E-autos}(b) below.

\begin{exa}\label{E-autos}
\begin{alphalist}
\item Let
$\F=\F_4=\{0,1,\alpha,\alpha^2\}$ and $n=3$. Then we compute
$x^n-1=\pi_1\pi_2\pi_3$ where $\pi_1=x+1,\ \pi_2=x+\alpha$ and $
\pi_3=x+\alpha^2$. In this case $s=1,\,r_1=3$, and $L_1=\F$. Thus
Corollary~\ref{autnumb} gives us $r_1!=6$ automorphisms, which are also
given in Example~\ref{E-Azs1}.
They all arise from pure permutations of the components.
\item Let $\F=\F_4$ as before and $n=5$. In this case
$$x^5-1=(x+1)(x^2+\alpha x+1)(x^2+\alpha^2x+1)$$ and we find
$s=2,\,r_1=1,\,r_2=2,\,L_1=\F,\,L_2\cong\F_8$.
Furthermore
\begin{equation}\label{e-K1-3}
  K_1=\F[x]/\ideal{x+1},\quad K_2=\F[x]/\ideal{x^2+\alpha x+1},\quad
  K_3=\F[x]/\ideal{x^2+\alpha^2x+1}.
\end{equation}
Corollary~\ref{autnumb} now says, that there are $ 1^1 2^21!\, 2!=8$
automorphisms. Once given the only nontrivial
$\F$-automorphism $\lambda$ of $L_2$ they can be listed systematically
according to Theorem~\ref{autos}. We want to present these
automorphisms with respect to the various descriptions of~$A$ as
in\eqnref{E-A},\eqnref{E-Adecomp}, and\eqnref{E-Ldecomp}. In order
to do so we first notice that $\lambda$ is given by the Frobenius homomorphism,
i.~e. $\lambda(a)=a^4$ for all $a\in L_2$.
Secondly, we need an $\F$-isomorphism between the
two fields $K_2$ and $K_3$. The list given below is
based on the isomorphism
\begin{equation}\label{e-psi}
  \Psi:\F[x]/\ideal{x^2+\alpha x+1}\longrightarrow
  \F[x]/\ideal{x^2+\alpha^2x+1},\quad\text{where }
  \Psi(x)= \alpha^2 x+1,\quad
\end{equation}
with inverse given by $\Psi^{-1}(x)=\alpha x+\alpha$. Going through
all the necessary isomorphisms one obtains the descriptions for the
automorphisms on~$A$ as given in the table below. In the first column
of the list we use the standard notation $(\rho_1,\rho_2,\rho_3)$ for
a permutation $\rho\in S_3$. In the second (resp.\ third) column the
image of $[1,x,x]$ (resp.~$x$) under the corresponding automorphism
is given. Recall that this fully determines the $\F$-automorphism.
For instance, the second column of the seventh row is obtained as
follows (in suggestive notation)
\begin{align*}
  \big((\id,\lambda,\id)\circ(1,3,2)\big)[1,\,x,\,x]&=
  (\id,\lambda,\id)[1,\,\Psi^{-1}(x),\,\Psi(x)]\\
  &=(\id,\lambda,\id)[1,\,\alpha x+\alpha,\,\alpha^2x+1]\\
  &=[1,\,(\alpha x+\alpha)^4~\mod(x^2+\alpha x+1),\,\alpha^2x+1].
\end{align*}
Hence this automorphism maps $[a,b,c]$ onto
$[a,\Psi^{-1}(c)^4,\Psi(b)]$.
The relation between the third and second column is given
by the Chinese Remainder Theorem, see\eqnref{E-rho}.
{\small
\[
\begin{array}{|c|c|c|}
\hline \F\times\F_8^2 &
\F[x]/{\scriptstyle\langle{x+1}\rangle}\times\F[x]/{\scriptstyle\langle{x^2+\alpha
x+1}\rangle}\times \F[x]/{\scriptstyle\langle{x^2+\alpha^2
x+1}\rangle} &
\F[x]/{\scriptstyle\langle{x^5-1}\rangle}\rule[-.5ex]{0cm}{3ex}\\
\hline\hline \!(\id,\id,\id)\circ(1,2,3)\!\!& [1,x,x] & x
\rule[-.5ex]{0cm}{3ex}\\
\hline \!(\id,\id,\lambda)\circ(1,2,3)\!\! &
[1,x,x^4]=[1,x,x+\alpha^2] & \alpha x^4+x^3+x^2+\alpha^2x
\rule[-.5ex]{0cm}{3ex}\\
\hline \!(\id,\lambda,\id)\circ(1,2,3)\!\! &
[1,x^4,x]=[1,x+\alpha,x] &\alpha^2x^4+x^3+x^2+\alpha x
\rule[-.5ex]{0cm}{3ex}\\
\hline \!(\id,\lambda,\lambda)\circ(1,2,3)\!\! &
[1,x+\alpha,x+\alpha^2] & x^4 \rule[-.5ex]{0cm}{3ex}\\
\hline\!(\id,\id,\id)\circ(1,3,2)\!\! & [1,\alpha x+\alpha,\alpha^2x+1] &
x^4+\alpha x^3+\alpha^2x^2+x \rule[-.5ex]{0cm}{3ex}\\ \hline
\!(\id,\id,\lambda)\circ(1,3,2)\!\! & [1,\alpha x+\alpha,(\alpha^2x+1)^4]=
            [1,\alpha x+\alpha,\alpha^2x+\alpha^2] & x^3
            \rule[-.5ex]{0cm}{3ex}
\rule[-.5ex]{0cm}{3ex}\\
\hline \!(\id,\lambda,\id)\circ(1,3,2)\!\! &
[1,(\alpha x+\alpha)^4,\alpha^2x+1]=
            [1,\alpha x+1,\alpha^2x+1] & x^2
\rule[-.5ex]{0cm}{3ex}\\
\hline \!(\id,\lambda,\lambda)\circ(1,3,2)\!\! &
[1,\alpha x+1,\alpha^2x+\alpha^2] & x^4+\alpha^2x^3+\alpha x^2+x
\rule[-.5ex]{0cm}{3ex}\\ \hline
\end{array}
\]
}
\end{alphalist}
\end{exa}

In the examples of the next two sections about the left ideals in $\Azs$, we will mainly use
a representation as displayed in the second column above.
Only thereafter, when dealing with cyclic codes, we will need
computations mod~$(x^n-1)$ as in the third column.

In the foregoing example~(b) the first four automorphisms do not
permute the components of~$A$.
In such a case there exist no non-trivial $\sigma$-CCC's as we will see in
part~(a) of the following result, which also can be regarded as an extension of
Proposition~\ref{P-trivialCCC}.
The if-part of this statement and part~(b)
can also be found in~\cite[Thm.~8 and Thm.~6]{Ro79}.

\begin{prop}\label{P-blockcode}
\begin{alphalist}
\item Let $\sigma \in \AutF(A)$ and $K^{(k)}$ be as in\eqnref{E-Kk}.
      Then every $\sigma$-CCC is a block code iff
      $\sigma(K^{(k)})=K^{(k)}$ for all $\ 1\le k\le r$.
\item Let $\cC$ be $\sigma$-CCC, $\cJ=\p (\cC )$
      be the corresponding ideal in $\cR$ and
      $\cJ_0:=\{c\in A\mid \exists\; g\in\cJ: g_0=c\}$, where $g_0$
      denotes the $z$-free term of~$g$.
      If $\sigma (\cJ_0 )=\cJ_0 $, then~$\cC$ is a block code.
\end{alphalist}
\end{prop}

It is possible to give a direct proof of the result at this point.
Since we don't need the proposition, it is most efficient to postpone the proof
to the end of Section~\ref{S-gencontr}.

The proposition demonstrates that an essential ingredient of a
nontrivial $\sigma$-CCC is the way of how~$\sigma$ properly
permutes the components of~$A$.
This in turn determines to a large extent the structure of the
algebra~$\cR=\Azs$.
To give an idea of this we mention without proofs the following facts
(which won't be used in the paper):
\begin{arabiclist}
\item Let $\bigcup_{j=1}^s Z_j$ be a partition of
      $\mathcal K:=\{K^{(1)},\ldots,K^{(r)}\}$ which is invariant under~$\sigma$,
      i.~e. $\sigma (Z_j)=Z_j$ for all $1\le j\le s$.
      Then~$\cR$ is a direct sum of subalgebras
      $\cR=\bigoplus_{j=1}^s \cR^{(Z_j)}$, where
      $\cR^{(Z_j)}=\sum_{K^{(i)} \in Z_j}\ve{i}\cR$.
\item Whenever~$Z_j$ contains exactly one field $K^{(i)}$,
      then~$\ve{i}\cR$ is a classical skew-polynomial domain.
\end{arabiclist}

\section{Generators for left ideals in
\protect{$\Azs$}}\label{S-idealgen} \setcounter{equation}{0}

As a first fundamental property we note that $\cR=\Azs$ inherits from $\F[z]^n$
the property 'left Noetherian' by means of the left $\F[z]$-isomorphism $\p $ from\eqnref{E-p}.
This is also a straightforward consequence
of results in Section~\ref{S-circ}, where $\cR$ will appear as the image
of $\F[z,x]$  under a left $\F[z]$-homomorphism
(see the discussion following Theorem~\ref{T-sMgsMh}).
In a similar way or by an anti-isomorphism as given below in Observation~\ref{O-anti}
one can see that $\cR$ is also right Noetherian. \\
The central theme in Piret's fundamental article~\cite{Pi76}
is the detailed construction of a generator polynomial for an
irreducible $\sigma$-CCC, resp.\ left ideal in $\cR$.
This is done for an
automorphism $\sigma$ which maps $x$ onto a power of $x$. The
constructions are displayed in terms of involved matrix manipulations.
But at the same time central arguments rely heavily on the
decomposition of $A$ into components as introduced in the foregoing
section. Maybe this is the reason why the small step in~\cite{Kae81}
for obtaining a single generator polynomial for {\em reducible\/} CCC's
is not done in~\cite{Pi76}. In this section we will first show by quite
different, rather short and purely algebraic arguments and for an
arbitrary automorphism $\sigma $ that any delay-free left ideal in
$\cR$ is in fact a  principal left ideal (Theorem~\ref{existence}).
This result is not constructive. The development of an algorithmic
procedure is postponed to the next section.

In~\cite{Pi76,Kae81} uniqueness of generator polynomials is not addressed. The key
to our uniqueness result in Theorem~\ref{T-unique} is a reduction procedure
which resembles the one in Groebner basis theory but which has to
take into account that $\Azs $ usually has many zero divisors and is
not commutative. It turns out that reduced generators are essentially
unique. At the same time reduced generators behave well for
explicitly writing down a generator matrix for the corresponding code
(see Section~\ref{S-gencontr}). They also will lead
directly to minimal generator matrices for CCC's. We conclude the section
with some information on right ideals which will be of later use, too.

In this section any isomorphic representation of~$A$ as a direct product of
fields with the corresponding unique set of pairwise orthogonal primitive
idempotents $\ve{1},\dots,\ve{r}$ will do.
A canonical way of displaying the algebra has been
described in\eqnref{E-xn-1} --\eqnref{E-Adecomp}.
However, in any case we obtain the fields (see also\eqnref{E-Kk} for the
canonical representation)
\[
  K^{(k)}:=\ve{k}A \text{ for } k=1,\ldots,r.
\]
The primitive idempotents will play a central role in the arguments of
this and the next section. Notice that $\sum_{k=1}^r\ve{k}$ is the
identity in~$A$, and thus also in~$\cR$, and therefore,
\begin{equation}\label{e-efk}
  f=\ve{1}f+\cdots+\ve{r}f \text{ for all }f\in\cR.
\end{equation}
Before we proceed let us introduce the following useful notation.

\begin{notation}\label{Nota1}
\begin{arabiclist}
\item For $f\in\cR$ and $k=1,\ldots,r$ put $f^{(k)}:=\ve{k}f$.
      We call $f^{(k)}$ the $k$-th component of~$f$. Furthermore, we call
      $T_f:=\{k\mid f^{(k)}\not=0\}$ the support of~$f$.
\item For a polynomial $f=\sum_{\nu=0}^d z^{\nu}f_{\nu}$, where
      $f_{\nu}\in A$ and $f_d\ne0$, we call $\deg_z f:=d$ the $z$-degree,
      $f_d$ the leading $z$-coefficient, and $f_0$ the $z$-free term
      of~$f$.
\item The left (resp.\ right) ideal in~$\cR$ generated by a set
      $M\subseteq\cR$ will be denoted by $\lideal{M}$
      (resp.\ $\rideal{M}$).
\end{arabiclist}
\end{notation}

From\eqnref{e-efk} we immediately obtain the following.

\begin{obs}\label{O-lideal}
Let $f_1,\ldots,f_t\in\cR$ and put $f_i^{(k)}:=\ve{k}f_i$ for
$i=1,\ldots,t$ and $k=1,\ldots,r$. Then
\[
  \lideal{f_1,\ldots,f_t}=
  \lideal{f_1^{(1)},\ldots,f_1^{(r)},\ldots\ldots,f_t^{(1)},\ldots,f_t^{(r)}}.
\]
\end{obs}

It is an elementary, but crucial fact that each automorphism
$\sigma\in\AutF(A)$ induces a permutation on the set of primitive
idempotents, i.~e.
\begin{equation}\label{e-perm}
  \{ \ve{1},\dots ,\ve{r}\}=
  \{ \sigma (\ve{1}),\dots ,\sigma (\ve{r})\}.
\end{equation}
This implies that for a given polynomial
$g=\sum_{\nu\geq0}z^{\nu}g_{\nu}\in\cR$ the $z$-coefficients of
the components
\begin{equation}\label{E-gk}
   \ve{k}g=\sum_{\nu=0}^dz^{\nu}\sigma^{\nu}(\ve{k})g_{\nu}
\end{equation}
are in general not in $K^{(k)}$ but rather move around according to the
permutation\eqnref{e-perm}.
In particular, for each $\nu\geq0$ and each $k\in\{1,\ldots,r\}$
there exists a unique $l\in\{1,\ldots,r\}$ such that
$\sigma^{\nu}(\ve{k})g_{\nu}\in K^{(l)}$, see also Example~\ref{E-norming} below.

The following lemma will be of frequent use.

\begin{lemma}\label{norming}
\begin{alphalist}
\item The element $x$ is a unit (i.e. invertible) in $A$ and
      $a\in A$ is a unit in~$A$ if and only if
      $\ve{k}a\ne 0$ for all $k\in\{1,\ldots,r\}$.
\item Let $g=\sum_{\nu=0}^dz^{\nu}g_{\nu}\in \cR$ and suppose that for some
      $\mu\in\{0,\dots ,d\}$ the coefficient $g_{\mu}$ is nonzero. Then there is
      a unit $a\in A$ such that
      $ag=\sum_{\nu=0}^d z^{\nu}(ag)_{\nu}=\sum_{\nu=0}^d z^{\nu}\sigma^{\nu}(a)g_{\nu}$,
      satisfies for all $1\le k\le r$
      $$\ve{k}g_{\mu}\ne 0\ \Longrightarrow\
       \ve{k} (ag)_{\mu} =\ve{k}.
      $$
\item Let $g\in\cR$ be a nonzero polynomial.
      Then there exists a unit $a\in A$, such that for all $1\le k\le r$
      \[
         k\in T_g\Longrightarrow
        \text{ the leading $z$-coefficient of $(ag)^{(k)}$ is a primitive idempotent.}
      \]
      We say that the polynomial $ag$ is normalized.
\item  Let $f,g\in\cR$ and $1\le k\le r$, then
      $$f\ve{k}g=0\ \Longrightarrow\ f\ve{k}=0 \text{ or } \ve{k}g=0.$$
\end{alphalist}
\end{lemma}

Notice that, since~$A$ is commutative, $(ag)^{(k)}=ag^{(k)}$ for
all $a\in A$ and $g\in\cR$.

\begin{proof}
(a) is obvious.
\\
(b) Since $K^{(k)}=\ve{k}A$ is a field, we can
find $b_k\in A$ such that
$\ve{k}b_kg_{\mu}=\ve{k}$ whenever $\ve{k}g_{\mu}\ne 0$.
Now
\[
   a:=\sigma^{-\mu}\Big(\sum_{k\in T_{g_{_{\scriptscriptstyle\mu}}}}\ve{k}b_k
      +\sum_{k\not\in T_{g_{_{\scriptscriptstyle\mu}}}}\ve{k}\Big).
\]
has the desired properties.
Invertibility follows from~(a).
\\
(c)
By the previous part we can find for each $k\in T_g$ units $a_k\in A$
such that $a_kg^{(k)}$ has a
primitive idempotent as leading $z$-coefficient. Let
$$a:=\sum_{k\in T_g}a_k\ve{k}+\sum_{k\not\in T_g}\ve{k}.$$
Then one easily verifies that~$a$ is a unit in~$A$ and
$ag^{(k)}=a_kg^{(k)}$ yields the desired property.
\\
(d)
If $f\ve{k}\not=0\not=\ve{k}g$, then the leading $z$-terms of $f\ve{k}$ and $\ve{k}g$
are of the form $z^{\nu}a\ve{k}\ne0$ and $\ve{k}bz^{\mu}\ne0$, respectively,
for some $a,\,b\in A$.
But then the leading $z$-term of $f\ve{k}g$ is
$z^{\nu}a\ve{k}bz^{\mu}$, which is nonzero by\eqnref{E-idemzero}.
\end{proof}

Note that part~(d) above extends\eqnref{E-idemzero}.

\begin{exa}\label{E-norming}
Let us consider the case $\F=\F_4$ and $n=5$.
The ring~$A$ and its automorphisms have been described in detail in
Example~\ref{E-autos}(b).
We now choose the automorphism~$\sigma$ given by $\sigma(x)=x^2$.
The effect of normalization is best visualized when representing the
elements in~$A$ as triples in $K_1\times K_2\times K_3$, where
the fields~$K_i$ are as in\eqnref{e-K1-3}, see also the list in
Example~\ref{E-autos}(b).
In this description we have
$\sigma([u,v,w])=[u,\Psi^{-1}(w)^4,\Psi(v)]$, where~$\Psi$ is as
in\eqnref{e-psi}.
The primitive idempotents
$\ve{1}=[1,0,0],\,\ve{2}=[0,1,0]$, and $\ve{3}=[0,0,1]$
satisfy $\sigma(\ve{1})=\ve{1},\,\sigma(\ve{2})=\ve{3}$, and
$\sigma(\ve{3})=\ve{2}$.
Consider now the element
\[
   g=[0,z(\alpha^2x+\alpha^2),\alpha x+1].
\]
Then one easily verifies that $\ve{1}g=\ve{2}g=0$ and $\ve{3}g=g$.
We want to normalize~$g$.
Since $(\alpha^2x+\alpha^2)^{-1}=x+\alpha$ in the field~$K_2$, we put
$a:=\sigma^{-1}([1,x+\alpha,1])$, see the proof of part~(b) above.
Since $\sigma^{-1}(x)=x^3$, or, in the current representation,
$\sigma^{-1}([u,v,w])=[u,\Psi^{-1}(w),\Psi(v)^4]$, we calculate
$a=[1,1,\alpha^2 x]$.
Now one checks that
\[
  ag=z[0,1,0]+[0,0,1]=[0,z,1].
\]
In this case normalization of the leading $z$-coefficient
led to a normalization of the $z$-free term, too.
\mbox{}\hfill$\Box$
\end{exa}

We can now proceed to our algebraic (generalized and completed) version
of Piret's result on ideal generators for $\sigma$-CCC's,
see~\cite[Thm.~3.10]{Pi76}.

\begin{theo}\label{existence}
Let $\cC$ be an $\F[z]$-submodule of $\F[z]^n$ and $\cJ=\p (\cC)$ its
image in $\cR$. Then the following properties are equivalent.
\begin{alphalist}
\item $\cC$ is $\sigma$-cyclic and delay-free.
\item $\cJ=\lideal{g}$ for some polynomial $g\in \cR$ satisfying
      $T_g=T_{g_{_{\scriptscriptstyle 0}}}$, precisely
      \begin{equation}\label{E-existence}
       \ve{k}g\ne 0 \Longleftrightarrow \ve{k}g_0 \ne 0
       \text{ for all } 1\le k\le r.
      \end{equation}
      Here $T_g$ denotes the support and $g_0$ the $z$-free term
      of~$g$, see~\ref{Nota1}.
\end{alphalist}
In particular, every delay-free left ideal of~$\cR$ is principal.
\end{theo}

\begin{proof}
For any polynomial $f\in\cR$ we will use the notation $f_0$ for its
$z$-free term.
\\
``(a)~$\Rightarrow$~(b)''
First of all, $\cJ$ is a left ideal by Observation~\ref{Piretalgebra}(b).
Thus it remains to show that~$\cJ$ has a principal generator
satisfying\eqnref{E-existence}.
For $k\in \{ 1,\dots ,r\} $ define $g^{(k)}:=0$ if
$\left (\ve{k}\cJ \right )_0:=\left\{f_0\mid f\in \ve{k}\cJ \right\} =\{\, 0\}$
and let otherwise
$g^{(k)}$ be a polynomial of minimal $z$-degree in $\left\{ f\in
\ve{k}\cJ\mid f_0\ne 0\right\}$. Multiplying by an appropriate
constant factor according to Lemma~\ref{norming}(b), we can assume
that $g^{(k)}_0=\ve{k}$ whenever $g^{(k)}\ne 0$.
Then for each $k\in\{1,\ldots,r\}$ either $g^{(k)}=0$ or
$g^{(k)}=\ve{k}+ \sum_{i=1}^{d_k}z^ig^{(k)}_i$
for some $d_k\ge 0$, $g_i^{(k)}\in A$, and $g^{(k)}_{d_k} \ne
0$. Put
\begin{equation}\label{gdef}
g:=\sum_{k=1}^rg^{(k)}.
\end{equation}
Obviously, $g^{(k)}=\ve{k}g$ and the notation matches
with~\ref{Nota1}(1).
By construction we have $\lideal{g}\subseteq\cJ$ as well as
property\eqnref{E-existence}.
Hence it remains to show that $\cJ\subseteq\lideal{g}$.
In order to do so, define the length of an arbitrary polynomial
$f=\sum_{i=i_0}^{i_0+d}z^if_i \in \cR$ with $f_{i_0}\ne 0 \ne f_{i_0+d}$
as $l(f):=d+1$  and put $l(0):=0$. Suppose
now that $\cJ \setminus \lideal{g} \ne \varnothing$ and let $f$ be a
polynomial of minimal length in $\cJ \setminus \lideal{g} $. We have
$f=z^{i_0}\overline{f}$ and $l(f)=l(\overline{f})$ for some
$\overline{f}\in\cR$ such that $\overline{f}_0\not=0$.
Delay-freeness of~$\cJ$, see Observation~\ref{O-attrib}, implies
$\overline{f}\in\cJ$, too.
Since $\overline{f}\notin \lideal{g}$ we can assume
without restriction $\overline{f}=f$, i.~e. $f_0\ne 0$. Now let
$f{'}:=f-f_0g$. Then $f{'}\in\cJ\backslash\lideal{g}$.
Moreover, we obtain for each $k\in T_g$ the identity
$$ \ve{k}f{'}=\ve{k}f-\ve{k}f_0g=\ve{k}f-f_0g^{(k)}.
$$
Since $g^{(k)}_0=\ve{k}$ we conclude $(\ve{k}f')_0=0$.
If $k\not\in T_g$, then $g^{(k)}=0$, which means that
$(\ve{k}\cJ)_0=\{0\}$, and thus $(\ve{k}f')_0=0$, too.
Hence $f'_0=0$ and, by the choice of the polynomials~$g^{(k)}$,
one has $\deg_z f'\leq\deg_z f$.
This together implies $l(f')<l(f)$, which contradicts the choice of
$f$. This proves $\cJ=\lideal{g}$.
\\
``(b)~$\Rightarrow$~(a)''
Let $\cJ=\lideal{g}$ be a principal left ideal and~$g$ satisfy\eqnref{E-existence}.
Since by Observation~\ref{Piretalgebra} $\cC=\v(\cJ)$ is $\sigma$-cyclic,
it remains to show that $\cJ$ is delay-free, see also Observation~\ref{O-attrib}.
In order to do so, we may assume by
Lemma~\ref{norming}(b) and\eqnref{E-existence} that for all
$1\le k\le r$ either $g^{(k)}_0=\ve{k}$ or $g^{(k)}=0$. Let now
$f=ug\in \cJ$ where $u=\sum_{\mu=0}^{\delta}z^{\mu}u_{\mu}$ and assume $f_0=0$.
Then $f=zf'$ for some $f'\in\cR$ and we have to show that $f'\in\cJ$.
From the equation
$$0=f_0=u_0g_0=u_0\sum_{k=1}^rg^{(k)}_0
=\sum_{k=1}^ru_0\ve{k}g_0^{(k)}=\sum_{k\in T_g}u^{(k)}_0$$
we get $u_0^{(k)}=0$ for all $k\in T_g$.
This in turn implies $u_0g=\sum_{k=1}^ru_0\ve{k}g=0$ and thus
$f=u{'}g$ for $u{'}:=u-u_0$. But
$u{'}=zu{'}{'}$ and we finally conclude that $f'=u''g\in \cJ$,
showing that $\cJ$ is delay-free.
\end{proof}

The next example shows that not all left ideals in $\cJ$ are
principal and that not every generator $g$ of a delay-free principal left
ideal fulfills\eqnref{E-existence}.

\begin{exa}
\begin{alphalist}
\item Let $\F=\F_4$ and  $n=3$ be as in Example~\ref{E-CCC1}. In Example~\ref{E-autos}~(a)
      we saw that up to an isomorphism $A=K_1\times K_2\times K_3$,
      where $K_i=\F[x]/\ideal{\pi_i}$.
      We choose the automorphism~$\sigma$ which corresponds to the permutation $(1,3,2)$.
      In the representation $A\cong \F[x]/\ideal{x^3-1}$ this
      corresponds to the automorphism which maps $x$ onto $\alpha^2x$.
      Now let $f_1=z[1,1,1],\,f_2=[0,1,0]$ and assume that
      $\lideal{f_1,f_2}=\lideal{g}$ for some $g\in\cR$.
      Then the $z$-free term of~$g$ is of the form $[0,a,0]$ for some
      $a\in A\backslash\{0\}$ and comparing $z$-coefficients in an equation
      $f_1=ug,\,u\in\cR$, leads to  a contradiction.
      Thus the left ideal $\lideal{f_1,f_2}$ is not principal.
      The same example works, mutatis mutandis, for any automorphism
      $\sigma\in\AutF(A)$ satisfying $\sigma(\ve{2})=\ve{3}$ and for
      any~$n$ and~$\F$ where, as usual, $\text{char}(\F)\nmid n$.
\item Let now $A$ be arbitrary and $\sigma\in\AutF(A)$ such that
      $\sigma (\ve{1})=\ve{2}$. Let $g=(z+1)\ve{2}$. Then
      $\sigma (\ve{2})\ne\ve{2}$ and the left ideal $\lideal{g}=\lideal{\ve{2}}$ is
      delay-free, but $\ve{1}g=z\ve{2}\ne 0$ and $\ve{1}g_0=0$.
      \mbox{}\hfill$\Box$
\end{alphalist}
\end{exa}

The proof of Theorem~\ref{existence} is not constructive as long as
there is no finite procedure to determine the minimal polynomials
$g^{(k)}\in \ve{k}\cJ$ starting from a finite generating family of~$\cJ$.
In the next section such a procedure will be developed. But before we go
into the computational issues we will investigate, as to what extent a
generator of a principal left ideal is unique.
The key to our uniqueness result is a reduction procedure based on a
monomial ordering which we introduce now.

\begin{defi}\label{MO}
\begin{alphalist}
\item For $\mu\ge 0$ and $1\le k\le r$ the polynomials $z^\mu\ve{k}$
      are called the (left-) monomials of~$\cR$.
\item Given two monomials $z^\mu\ve{k}, z^\nu\ve{l}$ we define
      $$ z^\mu\ve{k}< z^\nu\ve{l}:\Longleftrightarrow
      \mu<\nu \text{ or }
      \mu=\nu \text{ and } k<l.$$
\item Let $f=\sum_{\mu=0}^dz^{\mu}f_{\mu}\in\cR$ be nonzero and
      have the following component expansion
      \[
       f=\left( \ve{1}f_0+\cdots +\ve{r}f_0\right)+
       z\left(\ve{1}f_1+\cdots +\ve{r}f_1\right ) +
       \cdots + z^d\left(\ve{1}f_d+\cdots +\ve{r}f_d \right).
      \]
      Then the individual summands $z^{\mu}\ve{k}f_{\mu}$,
      as far as nonzero, are called the terms of $f$.
      The (left-) leading monomial, denoted by
      $\LM{f}$, is the largest monomial $z^{\mu}\ve{k}$
      (with respect to~$<$) such that $\ve{k}f_{\mu}\ne 0$.
      The associated term is called the leading term of~$f$.
\end{alphalist}
\end{defi}

Observe that in the canonical representation of~$A$ as given
in\eqnref{E-Adecomp},\eqnref{E-canIdem} the monomials are of the form
\[
    z^{\mu}\ve{k}=[0,\ldots,0,z^{\mu},0,\ldots,0]
    \text{ where $z^{\mu}$ is at the $k$-th position.}
\]
In the context of ordinary Groebner basis theory (i.e. commutative
and no zero-divisors) one would like to call such an ordering a
TOP-monomial ordering (Term Over Position) and in fact one
readily verifies that~(b) defines a well-ordering on the set of all
monomials which respects left-multiplication by monomials as far as
the result is nonzero. As will soon become clear our results actually
will not depend on the way the components are ordered in the
representation of~$A$.

Since we won't make use of any right monomials we will call
left monomials simply monomials. The following rules will be very useful.

\begin{lemma}\label{L-monorule}
\begin{alphalist}
\item For all $a\in A$ and all possible $\mu,\,\nu,\,k,\,l$ one has:
      \\
      $z^{\mu}\ve{k}\text{ is a right divisor of }z^{\nu}\ve{l}a
         \Longleftrightarrow \mu\le \nu \text{ and } k=l$.
\item Let $g,\,g'\in\cR$ and $k,\,l\in\{1,\ldots,r\}$ such that
      $\ve{k}g\ne0\ne\ve{l}g'$ and $\LM{\ve{k}g}=z^{\alpha}\LM{\ve{l}g'}$
      for some $\alpha\in\N_0$.
      Then $\ve{l}=\sigma^{\alpha}(\ve{k})$.
\end{alphalist}
\end{lemma}

\begin{proof}
Part~(a) is obvious. As for~(b), let $\deg_z(\ve{k}g)=d$ and
$\deg_z(\ve{l}g')=d'$.
Then $\LM{\ve{k}g}=z^d\sigma^d(\ve{k})$ and
$\LM{\ve{l}g'}=z^{d'}\sigma^{d'}(\ve{l})$.
Now the assumption implies $d=\alpha+d'$ and
$\sigma^{\alpha+d'}(\ve{k})=\sigma^d(\ve{k})=\sigma^{d'}(\ve{l})$, from
which the assertion follows.
\end{proof}

Now we turn to the notion of reducedness.

\begin{defi}\label{reduce}
\begin{alphalist}
\item Let $f_1,\dots, f_s$ be any family of polynomials from $\cR$.
      The family is called
      (left-) reduced if for all $1\le k,l\le s$ such that $k\ne l$ and
      $f_k\ne 0 \ne f_l$
      no nonzero term of $f_k$ is right divisible by $\LM{f_l}$.
\item A single polynomial $g\in \cR$ is called (left-) reduced if
      the family $\ve{1}g,\dots ,\ve{r}g$ is left-reduced.
\end{alphalist}
\end{defi}

Again, we will usually skip the qualifier `left'.
Note that a reduced family might contain one or more zero polynomials, but, of
course, no other polynomial appears more than once.

The following result describes the basic reduction process which will
lead us to unique ideal generators.
They will later on turn out to have further nice properties.
For the process we will need so-called 'elementary operations' on a family
$f_1,\ldots,f_s$, by which we mean the replacement of some~$f_k$ by
\begin{equation}\label{E-elop}
f'_k:=f_k-z^{\mu}af_l, \text{ for any } l\ne k,\mu\in\N_0 \text{ and } a\in A.
\end{equation}

\begin{prop}\label{redu}
Any finite family $f_1,\,\dots,\,f_s$ from $\cR$ can be
transformed by finitely many elementary operations into a
reduced family $g_1, \dots ,g_s$ such that
for the respective left ideals one has
$$\lideal{f_1,\dots ,f_s}=\lideal{g_1, \dots ,g_s}.$$
\end{prop}

\begin{proof}
First of all, it is clear that elementary operations leave the
corresponding left ideal invariant.
As for the reduction assume now that the leading term of some~$f_k$ is given by $z^{\nu}b$
and is right divisible by $\LM{f_l}$ for some $l\ne k$, say
\begin{equation}\label{redu0}
   z^{\nu}b=z^{\mu}\hat{a}\LM{f_l}\text{ for some }\hat{a}\in A,\,\mu\in\N_0.
\end{equation}
Define
\begin{equation}\label{redu1}
  f'_k:=f_k-z^{\mu}af_l,
\end{equation}
where $a\in A$ is such that\eqnref{redu0} holds true when we replace~$\hat{a}$ by~$a$ and
$\LM{f_l}$ by the leading term of~$f_l$.
Then either $f'_k=0$ or $\LM{f'_k}<\LM{f_k}$.
Observe also that $\deg_z(f'_k)\leq\deg_z(f_k)$ and equality is possible.
Proceed now with the family $f_1,\ldots,f'_k,\ldots,f_s$.
Since~$<$ is a well-ordering, we get after finitely many steps a family
$\hat{f}_1,\ldots,\hat{f}_s$, where no leading term is right divisible by
any other.
\\
As a second and final step we now autoreduce the family
$\hat{f}_1,\ldots,\hat{f}_s$.
Assuming that a term of $\hat{f}_k$, say $z^{\nu}b$, is right divisible by
$\LM{f_l}$ for some $l\ne k$, we proceed as in\eqnref{redu0}
and\eqnref{redu1}.
Since these operations do not affect the higher terms of $\hat{f}_k$ we
arrive after finitely many steps at the desired family.
\end{proof}

The following case of the reduction step will be of specific importance.

\begin{obs}\label{O-redu}
If in\eqnref{redu0} and\eqnref{redu1} $f_k\in\ve{k}\cR$ and
$f_l\in\ve{l}\cR$, then $f'_k\in\ve{k}\cR$, too.
This follows from the fact that $z^{\nu}b$ is a term of~$f_k=\ve{k}f_k$ and hence
$\ve{k}z^{\nu}=z^{\mu}\LM{f_l}$ by\eqnref{redu1}. Using
Lemma~\ref{L-monorule}(b), this implies $\ve{l}=\sigma^{\mu}(\ve{k})$ and
as a consequence
$f'_k=f_k-z^{\mu}af_l=\ve{k}f_k-z^{\mu}a\ve{l}f_l=\ve{k}(f_k-z^{\mu}af_l)\in\ve{k}\cR$.
\end{obs}

\begin{exa}\label{E-reduction}
Let us again consider the case $\F=\F_4$ and $n=5$ as
described in Example~\ref{E-autos}(b).
Just like in Example~\ref{E-norming} we will represent the elements as triples
in $K_1\times K_2\times K_3$, where
the fields~$K_i$ are as in\eqnref{e-K1-3}.
We now choose the automorphism~$\sigma$ given by
$\sigma([u,v,w])=[u,\Psi^{-1}(w)^4,\Psi(v)]$, where~$\Psi$ is as
in\eqnref{e-psi}, see the sixth line of the list in~\ref{E-autos}(b).
The primitive idempotents
satisfy $\sigma(\ve{1})=\ve{1},\,\sigma(\ve{2})=\ve{3}$, and
$\sigma(\ve{3})=\ve{2}$.
Consider the family $f_1,\ldots,f_6$, where
\begin{align*}
  f_1&=z\ve{1}=z[1,0,0],\\
  f_2&=z\ve{3}(\alpha x+\alpha)+\ve{2}\alpha^2x
     =z[0,0,\alpha x+\alpha]+[0,\alpha^2x,0],\\
  f_3&=z\ve{1}+\ve{1}=z[1,0,0]+[1,0,0],\\
  f_4&=z^2\ve{1}\alpha+z\ve{1}\alpha^2+\ve{1}=z^2[\alpha,0,0]+z[\alpha^2,0,0]+[1,0,0],\\
  f_5&=z\ve{3}(\alpha x+1)+\ve{2}(\alpha^2x+\alpha^2)=z[0,0,\alpha x+1]+[0,\alpha^2 x+\alpha^2,0],\\
  f_6&=z^2\ve{3}(\alpha^2x+\alpha^2)+z\ve{2}x=z^2[0,0,\alpha^2x+\alpha^2]+z[0,x,0].
\end{align*}
Note that in each case the first term is the leading term.
The family is not reduced and we perform the following steps.
\begin{arabiclist}
\item $f'_3:=f_3-f_1=\ve{1}$.
\item $f'_1:=f_1-zf'_3=0$.
\item $f'_4:=f_4-z^2\alpha f'_3=z\ve{1}\alpha^2+\ve{1}$.
\item $f''_4:=f'_4-z\alpha^2f'_3=\ve{1}$.
\item $f'''_4:=f''_4-f'_3=0$.
\item $f'_5:=f_5-af_2$, where $a\in A$ is such that
      $z\ve{3}(\alpha x+1)=az\ve{3}(\alpha x+\alpha)$.
      Hence $a=\sigma^{-1}[0,0,c]$, where $c=(\alpha x+1)(\alpha
      x+\alpha)^{-1}=\alpha^2x\in K_3$, and we get
      $a=[0,\Psi^{-1}(\alpha^2x)^4,0]=[0,x+\alpha^2,0]=\ve{2}(x+\alpha^2)$.
      Then we compute
      $f'_5=\ve{2}(\alpha^2x+\alpha^2)-\ve{2}(x+\alpha^2)\ve{2}\alpha^2x=0$.
\item $f'_6:=f_6-z\alpha f_2=0$.
\end{arabiclist}
Now the family
\[
  \hat{f}_1=0,\quad \hat{f}_2=f_2=z\ve{3}(\alpha x+\alpha)+\ve{2}\alpha^2x,\quad
  \hat{f}_3=\ve{1},\quad \hat{f}_4=0,\quad \hat{f}_5=0,\quad \hat{f}_6=0
\]
is reduced. We know that $\lideal{f_1,\ldots,f_6}=\lideal{\hat{f}_2,\,\hat{f}_3}$.
Applying Lemma~\ref{norming}(c), we can even normalize the generators and
obtain (after changing the ordering and omitting zero polynomials)
\[
 g_1:=\ve{1},\quad
 g_2:=\sigma^{-1}[1,1,(\alpha x+\alpha)^{-1}]\hat{f}_2=[1,\alpha
 x+\alpha^2,1]\hat{f}_2=z\ve{3}+\ve{2}.
\]
Since $g_k\in\ve{k}\cR$ for $k=1,\,2$ we know from Observation~\ref{O-lideal}
that $\lideal{f_1,\ldots,f_6}=\lideal{g}$, where
$g:=g_1+g_2=z[0,0,1]+[1,1,0]\in\cR$. Thus we have found a reduced and normalized
generator of the left ideal generated by $f_1,\ldots,f_6$.
\mbox{}\hfill$\Box$
\end{exa}

On first sight, the example appears somewhat specific in the sense that all given
generator polynomials are components, precisely
$f_1,\,f_3,\,f_4\in\ve{1}\cR$, $f_2,\,f_5\in\ve{2}\cR$, and $f_6\in\ve{3}\cR$.
However, by virtue of Observation~\ref{O-lideal} each ideal has a generating set
consisting of components only.

\begin{cor}\label{C-redu}
\begin{alphalist}
\item For every $f\in\cR$ there exists a unit $u\in\cR$, i.~e.
      $u\overline{u}=\overline{u}u=1$ for some $\overline{u}\in\cR$, such that
      the polynomial $uf$ is reduced.
      In particular, every principal left ideal has a reduced generator.
\item Every delay-free principal left ideal has a reduced generator
      satisfying property\eqnref{E-existence}.
\end{alphalist}
\end{cor}

\begin{proof}
(a) Only the first statement needs to be proven.
Define
\begin{equation}\label{e-fk}
    f_k:=\ve{k}f=f^{(k)}\text{ for }1\le k\le r.
\end{equation}
Then $f_k\in\ve{k}\cR$ and by definition~$f$ is reduced if and only if the
family $f_1,\ldots,f_r$ is reduced.
In order to prove the corollary we will analyze the effect of the
reduction process on the polynomial~$f$.
It suffices to consider a single reduction step as in\eqnref{redu0}
and\eqnref{redu1}, the result of which is the family
$f_1,\ldots,f'_k,\ldots,f_r$.
We will prove that
\begin{romanlist}
\item $f'_k\in\ve{k}\cR$,
\item $u:=1-z^{\mu}a\ve{l}$ is a unit in $\cR$,
\item $f':=f'_k+\sum_{j\ne k}f_j$ satisfies $f'=uf$.
\end{romanlist}
Part~(i) is in Observation~\ref{O-redu}. As a consequence, one
also has
\begin{equation}\label{e-zmuaf}
  z^{\mu}af_l\in\ve{k}\cR.
\end{equation}
As for~(ii), one easily derives from\eqnref{e-zmuaf}
that $u\overline{u}=\overline{u}u=1\text{ for }\overline{u}:=1+z^{\mu}a\ve{l}$.
Finally,~(iii) is established once we have shown that $\ve{j}f'=\ve{j}uf$ for all
$j=1,\ldots,r$. Using again\eqnref{e-zmuaf} and the orthogonality of the idempotents
we obtain for $j=k$ the identity
$\ve{k}uf=f_k-\ve{k}z^{\mu}af_l=f'_k=\ve{k}f'$ while for $j\ne k$
we have $\ve{j}uf=f_j-\ve{j}z^{\mu}af_l=f_j=\ve{j}f'$.
This completes the proof of~(a).
\\
(b) By Theorem~\ref{existence} we may assume that $\cJ=\lideal{f}$, where $f\in\cR$
satisfies\eqnref{E-existence}.
Again, it suffices to show that a single reduction step\eqnref{redu1}  respects this
property. But this is clear since ${f'}^{(j)}=f^{(j)}$ for all $j\not=k$
and\eqnref{e-zmuaf} shows that\eqnref{redu1}
occurs only for $\mu>0$ and in this case ${f'}^{(k)}_0=f^{(k)}_0$.
\end{proof}

In the proof we made use of the monomial ordering for the case where the
family consists of the components $f^{(1)},\ldots, f^{(r)}$ of a single
polynomial~$f\in\cR$.
In this case the leading term is uniquely determined
by the $z$-degree only and the arbitrarily prescribed ordering of the
idempotents $\ve{1},\ldots,\ve{r}$ has no effect on the reducedness.
It simply determines the ordering of the family $f^{(1)},\ldots,f^{(r)}$.

One should notice that in~(ii) of the proof above we encounter one of the many
units of the ring~$\cR$ which are not constant polynomials.

The following basic properties of reduced families will be of
essential use in the sequel.

\begin{lemma}\label{rules}
\begin{alphalist}
\item Let $g=g^{(k)}\in\ve{k}\cR $ and $u\in\cR$ such that $ug\ne 0$. Then
      $$\LM{ug}=z^{\alpha}\LM{g} \text{ for some }
      \alpha \ge 0.  $$
\item Let $G$ be a finite reduced subset of~$\cR$ such that
      $G=\bigcup_{k=1}^r G^{(k)}$, where $G^{(k)}:=\ve{k}G$; in other words,
      any element of $G$ is contained in one of the sets
      $\ve{k}\cR,\,k=1,\ldots,r$.
      Let furthermore $f\in\cR$. Then
      \[
         f\in\lideal{G}\Longrightarrow\LM{f}=z^{\alpha}\LM{g}
         \text{ for some $g\in G$ and some }\alpha\in\N_0.
      \]
\item Let $g\in\cR$ be a nonzero reduced polynomial and
      $u_1,\dots ,u_r\in\cR$. Then
      \[
          \sum_{k=1}^ru_kg^{(k)}=0\Longrightarrow u_kg^{(k)}=0\text{ for all }
          k=1,\ldots,r.
      \]
\end{alphalist}
\end{lemma}

\begin{proof}
(a)
Let $u=\sum_{\nu=0}^{\delta}z^{\nu}u_{\nu}$ and
$g=\ve{k}g=\sum_{\mu=0}^dz^{\mu}\sigma^{\mu}(\ve{k})g_{\mu}$, where
$\sigma^d(\ve{k})g_d\ne0$.
Then $ug=\sum_{\nu=0}^{\delta}z^{\nu}u_{\nu}g$ and if for some $\nu$
we find $u_{\nu}g\ne 0$, then also $u_{\nu}\ve{k}\ne 0$ and thus
\[
    \LM{u_{\nu}g}=\text{LM}\Big(\sum_{\mu=0}^d z^{\mu}\sigma^{\mu}(u_{\nu}\ve{k})g_{\mu}\Big)
    =\LM{g}
\]
since $\sigma^d(u_{\nu}\ve{k})\ne 0$ and by\eqnref{E-idemzero} also
$\sigma^d(u_{\nu}\ve{k})g_d\ne 0$. In order to find the leading monomial
of~$ug$ we thus only have to pick the maximal power $z^{\alpha}$ such that
$u_{\alpha}g\ne 0$ and then
$$\LM{ug}=\LM{z^{\alpha}u_{\alpha}g}=z^{\alpha}\LM{g}.$$
(b)
Suppose $G^{(k)}=\{g_1^{(k)},\dots ,g^{(k)}_{m_k}\}$, where $m_k=0$
if
$G^{(k)}=\varnothing$ and let
$$f=\sum_{k=1}^r\sum_{j=1}^{m_k}u_{kj}g_j^{(k)}
=\sum_{k=1}^r\sum_{j=1}^{m_k}
\underbrace{\left(u_{kj}\ve{k}\right)g_j^{(k)}}_{=:f_{kj}}
\text{ for some }u_{kj}\in\cR.$$
By part~(a) we have
$$\LM{f_{kj}}=z^{\beta_{kj}}\LM{g_j^{(k)}}
\text{ for some }\beta_{kj}\in\N_0\text{ if }f_{kj}\ne0.
$$
Consider now the leading monomials of the polynomials~$f_{kj}$ of maximal
$z$-degree. By reducedness of~$G$ these monomials are all different and this
proves the assertion.\\
(c)
Let $\sum_{k=1}^ru_kg^{(k)}=0$. By part (a) we know
already that $\LM{u_kg^{(k)}}=z^{\alpha_k}\LM{g^{(k)}}$ for some $\alpha_k\ge0$
whenever $u_kg^{(k)}\ne 0$. If there are nonzero products
$u_kg^{(k)}$ at all, then there must be some cancellation of the
maximal leading monomials which contradicts reducedness.
\end{proof}

We can now apply these techniques in order to obtain uniqueness
of generators of left ideals if we also assume normalization in the
sense of Lemma~\ref{norming}(c).

\begin{theo}\label{T-unique}
\begin{alphalist}
\item Every left ideal in~$\cR$ has a unique finite left-reduced generating
      family, such that each element is nonzero, normalized, and contained in one of the
      components $\ve{1}\cR,\,\ldots,$ $\ve{r}\cR$.
\item Every principal left ideal in~$\cR$ has a unique left-reduced and
      normalized generator.
\end{alphalist}
\end{theo}

\begin{proof}
Part~(b) is a consequence of~(a) and Corollary \ref{C-redu}.
\\
As for~(a) notice that, by virtue of Observation~\ref{O-lideal}, each left
ideal has a generating family consisting only of polynomials in the components
$\ve{1}\cR,\,\ldots,\,\ve{r}\cR$.
Using Observation~\ref{O-redu} this property is preserved when reducing
the family. Normalizing each element then proves the existence of the
desired generating family.
As for uniqueness, assume
\[
  \cJ=\lideal{g_1^{(k_1)},\ldots,g_s^{(k_s)}}=
  \lideal{{g'_1}^{(l_1)},\ldots,{g'_t}^{(l_t)}}
\]
where $g_i^{(k_i)}\in\ve{k_i}\cR$ and
${g'_i}^{(l_i)}\in\ve{l_i}\cR$, where all polynomials are normalized, and both families
are reduced.
Let $i\in\{1,\ldots,s\}$.
By Lemmata~\ref{rules}(b),~\ref{L-monorule}(b) one has
\begin{align*}
  \LM{g_i^{(k_i)}}&=z^{\alpha}\LM{{g'_j}^{(l_j)}}\text{ for some }
  \alpha\in\N_0\text{ and $j$ such that }
        \ve{l_j}=\sigma^{\alpha}(\ve{k_i}),\\[.5ex]
  \LM{{g'_j}^{(l_j)}}&=z^{\beta}\LM{g_m^{(k_m)}}\text{ for some }
  \beta\in\N_0\text{ and $m$ such that }\ve{k_m}=\sigma^{\beta}(\ve{l_j}).
\end{align*}
Combining these two equations we first obtain
$\LM{g_i^{(k_i)}}=z^{\alpha+\beta}\LM{g_m^{(k_m)}}$.
Then by reducedness of the family
$g_1^{(k_1)},\ldots,g_s^{(k_s)}$ we conclude successively
$\alpha=\beta=0 $, $i=m$, $k_i=k_m=l_j$.
Hence each leading monomial of the family $g_1^{(k_1)},\ldots,g_s^{(k_s)}$
occurs as leading monomial of the other family.
Symmetry and the fact that the leading monomials of a reduced family are
pairwise different, shows that $s=t$ and, after reordering,
\[
  k_i=l_i\;\text{ and }\; \LM{g_i^{(k_i)}}=\LM{{g'_i}^{(k_i)}}
  \text{ for all }i=1,\ldots,s.
\]
Suppose now that for some $i$ we have $g_i^{(k_i)}\ne {g'_i}^{(k_i)}$, or
equivalently $f:=g_i^{(k_i)}-{g'_i}^{(k_i)}\ne 0$. By normalization the
leading terms of $g_i^{(k_i)}$ and ${g'_i}^{(k_i)}$ are equal
and thus cancel. Therefore any term of $f$
comes from a non leading term
of $g_i^{(k_i)}$ or ${g'_i}^{(k_i)}$ or is a difference of such terms.
Since $f\in\cJ $, by Lemma~\ref{rules}(b) the leading term of $f$ must be
right divisible by some $\LM{g_j^{(k_j)}},\ j\ne i $. This contradicts
reducedness. Thus both reduced and normalized families must coincide.
\end{proof}

In the same way as $\sigma$-CCC's are linked to left ideals, their
duals will turn out to be linked to certain right ideals in $\cR$.
Our results on left ideals can be translated to right ideals by
means of the following anti-isomorphism.

\begin{obs}\label{O-anti}
For any $\sigma\in\AutF(A)$ the map \
$\widetilde{\text{ } }\ :\ \Azs\longrightarrow A[z;\sigma^{-1}]$
defined by
\begin{equation}
     g=\sum_{\nu\ge 0}z^{\nu}g_{\nu} \longmapsto \widetilde{g}
     :=\sum_{\nu\ge 0}g_{\nu}z^{\nu}=
     \sum_{\nu\ge 0}z^{\nu}\sigma^{-\nu}(g_{\nu})
\end{equation}
is an $\F$-algebra anti-isomorphism.
\end{obs}

Theorem~\ref{existence} immediately implies

\begin{cor}\label{c-anti}
Any delay-free right ideal $\cJ$ in $\Azs$ is a principal right ideal.
\end{cor}

In the next section the results will be
complemented by a computational procedure, which checks whether a
finitely generated left ideal is delay-free
or principal and, if so, computes the
unique generator polynomial.

\section{On the computation of principal generators of left ideals }\label{gencomp}
\setcounter{equation}{0}

While establishing uniqueness of a generator polynomial has been
(typically) somewhat more cumbersome, the computation --- starting from a
finite set of generators of a delay-free left ideal~$\cJ$ ---
can be achieved by a rather straightforward and systematic procedure.
Remembering the proof of Theorem~\ref{existence} it will be sufficient to
compute minimal $z$-degree polynomials with nonzero
constant term in each component $\ve{k}\cJ,\,1\le k\le r$, in order to obtain a
single generating polynomial.
Thereafter, reduction and normalization will lead to uniqueness according to
Theorem~\ref{T-unique} and Theorem~\ref{existence}.
As we will show in Theorem~\ref{T-algG} below, we obtain such minimal polynomials if we pick
any finite set of generators of the ideal, decompose it into its components, and
apply the reduction procedure to the family of components.
Furthermore, the algorithm even provides a test whether the ideal under
consideration is principal and/or delay-free.
The details are as follows.

Let $f_1,\dots ,f_s\in\cR=\Azs$ be any finite family and define
\begin{equation}\label{e-Fk}
   F^{(k)}:=\{f_1^{(k)},\ldots,f_s^{(k)}\}\text{ for }k=1,\ldots,r\text{ and }
   F:=\bigcup_{k=1}^r F^{(k)}.
\end{equation}
By Observation~\ref{O-lideal}
\begin{equation}\label{e-JF}
    \cJ:=\lideal{f_1,\ldots,f_s}=\lideal{F}.
\end{equation}
Note that some of the sets $F^{(k)}\subseteq \ve{k}\cJ$ may just
contain the zero polynomial but typically they also contain nonzero
polynomials.
It is quite surprising that just reducing the set~$F$ leads us, after
normalization, to the unique reduced and normalized generator polynomial.
The important observation here is, that when reducing
some polynomial $f_i^{(k)}\in F^{(k)}$ by some other $f_j^{(l)}\in F^{(l)}$,
the result $f_i^{(k)}-z^taf_j^{(l)}$ is necessarily again in
$F^{(k)}\subseteq \ve{k}\cJ$.
Therefore, the reduction process respects the partition
$F=\bigcup_{k=1}^r F^{(k)}$, and only the contents of
the individual sets $F^{(k)}$ changes.
The following theorem describes what can be obtained by reducing $F$.

\begin{theo}\label{T-algG}
Let~$F$ and~$\cJ$ be as in\eqnref{e-Fk} and\eqnref{e-JF}.
Furthermore, let~$F$ be transformed via finitely many elementary operations into the
reduced family~$G$. Define $G^{(k)}:=\ve{k}G$ for $k=1,\ldots,r$ and let
\[
    T:=\big\{k\in\{1,\ldots,r\}\,\big|\, G^{(k)}\ne\{0\}\big\}.
\]
Then
\begin{alphalist}
\item $G=\bigcup_{k=1}^r G^{(k)}$ and $\cJ=\lideal{G}$.
\item $\cJ$ is principal if and only if for each $k\in T$ the set $G^{(k)}$
      contains exactly one nonzero polynomial. Furthermore, if~$\cJ$ is
      principal then $\cJ=\lideal{g}$ where $g=\sum_{k\in T}g^{(k)}$
      and~$g^{(k)}$ is the unique nonzero polynomial in~$G^{(k)}$.
      In particular, the polynomial~$g$ is reduced.
\item $\cJ$ is delay-free if and only if~$\cJ$ is principal and the
      polynomial~$g$ of part~(b) satisfies\eqnref{E-existence}.
\item Let $\cJ$ be delay-free.
      Then for each $k\in T$ one has $\deg_z g^{(k)}\leq \deg_z f$ for all
      $f\in\ve{k}\cJ$ with nonzero $z$-free term.
\end{alphalist}
\end{theo}

\begin{proof}
(a) follows from Observations~\ref{O-redu} and~\ref{O-lideal}.
\\
(b) ``$\Leftarrow$'' and the second statement are consequences of Observation~\ref{O-lideal}.
``$\Rightarrow$'' Note that if $\cJ=\lideal{f}$, then
according to Observation~\ref{O-redu}, the
reduction of the
family $f^{(1)},\ldots,f^{(r)}$ leads to at most one
polynomial in each
component $\ve{k}\cR$, and therefore the assertion is a consequence of the
uniqueness in Theorem~\ref{T-unique}.
\\
(c) ``$\Leftarrow$'' is in Theorem~\ref{existence} while
``$\Rightarrow$'' is a combination of Corollary~\ref{C-redu}(b) and
Theorem~\ref{T-unique}.
\\
(d)
Suppose that for some $k\in T$ there exists a polynomial
$f\in\ve{k}\cJ$ satisfying $\text{deg}_zf<\text{deg}_zg^{(k)}$ and
having nonzero $z$-free term.
Then there is a constant $c\in A$ and a polynomial $\overline{f}\in\cR$ such that
$g^{(k)}-cf=z\overline{f}$.
By delay-freeness we have $\overline{f}\in\cJ$ and Lemma~\ref{rules}(b) implies
\begin{equation}\label{propG1}
   \LM{g^{(k)}}=z\LM{\overline{f}}=z^{1+\alpha}\LM{g^{(l)}}\text{ for some }
   l\in\{1,\ldots,r\}\text{ and }\alpha\geq0,
\end{equation}
contradicting reducedness of~$G$.
\end{proof}

Based on the forgoing proposition we have the following simple algorithmic
procedure for the computation of the unique reduced and normalized generator~$g$ of a
given delay-free ideal~$\cJ\subseteq\cR$.

\begin{algorithm}\label{Algo}
{\bf Input:}\quad A finite set $f_1,\dots , f_s$ of generators of the left ideal~$\cJ$.
\begin{algolist}
\item For all $1\le k\le r$ calculate $\ve{k}f_l,\,1\le l\le s$ and form the
      sets $F^{(k)}.$
\item Reduce the set $F=\bigcup_{k=1}^r F^{(k)}$ to obtain the reduced sets~$G^{(k)}$
      and~$G$.
\item Evaluation of results:
      \begin{caselist}
      \item If~$G^{(k)}$ contains more than one nonzero polynomial for some $k=1,\ldots,r$,
            then~$\cJ$ is not principal and thus not delay-free.
      \item If each set~$G^{(k)}$ contains at most one nonzero polynomial, denoted
            by~$g^{(k)}$, put $g:=\sum g^{(k)}$ and normalize~$g$ according to
            Lemma~\ref{norming}(c). Then $\cJ=\lideal{g}$ is principal
            and $g$ is its unique reduced and normalized generator.
            Furthermore,~$\cJ$ is delay-free if and only if~$g$
            satisfies\eqnref{E-existence}.
     \end{caselist}
\end{algolist}
\end{algorithm}

We close this section by an example.

\begin{exa}\label{E-comp}
Consider again the situation of Example~\ref{E-reduction} with the automorphism
given therein.
Furthermore, let $\cJ=\lideal{h_1,h_2,h_3}$, where
\begin{align*}
   h_1&:=z\ve{1}+\ve{2}\alpha^2x+z\ve{3}(\alpha x+\alpha),\\
   h_2&:=z\ve{1}+\ve{1}+\ve{2}\alpha^2x+z\ve{3}(\alpha x+\alpha),\\
   h_3&:=z^2\ve{1}\alpha\!+\!z^2\ve{3}(\alpha^2 x\!+\!\alpha^2)\!+\!z\ve{1}\alpha^2
         \!+\!z\ve{2}x\!+\!z\ve{3}(\alpha x\!+\!1)
         \!+\!\ve{1}\!+\!\ve{2}(\alpha^2 x\!+\!\alpha^2).
\end{align*}
As a first step we have to compute the components of these polynomials.
They are, not counting the zero components, just given by the polynomials in
Example~\ref{E-reduction}, precisely
\begin{align*}
   &\{h_1^{(1)},\,h_2^{(1)},\,h_3^{(1)}\}=\{f_1,\,f_3,\,f_4\},\ \;
   \{h_1^{(2)},\,h_2^{(2)},\,h_3^{(2)}\}=\{f_2,\,f_2,\,f_5\},\\[.5ex]
   &\{h_1^{(3)},\,h_2^{(3)},\,h_3^{(3)}\}=\{0,\,0,\,f_6\}.
\end{align*}
Thus, $\cJ=\lideal{f_1,\ldots,f_6}=\lideal{h}$, where the reduced and
normalized generator
\[
     h=z\ve{3}+\ve{1}+\ve{2}
\]
has already been calculated in Example~\ref{E-reduction}.
\mbox{}\hfill$\Box$
\end{exa}

\section{\protect{$\sigma$}-circulant matrices}\label{S-circ}
\setcounter{equation}{0}
While in the last sections we have concentrated on $\sigma$-CCC's
as left ideals in $\Azs$ we now focus on the description of
these codes as submodules of $\F[z]^n$.
More precisely, we introduce  $\sigma$-circulant matrices as
a counterpart of a generator polynomial of a principal left ideal.
These matrices show close resemblance  to classical circulants
which are common in the theory of cyclic block codes.
As a guideline through this section we first recall some basic facts
about classical circulant matrices over finite fields.
Many of these properties can then be generalized
appropriately  to $\sigma$-circulants.
The consequences for $\sigma$-CCC's will then be discussed in the next section.

Throughout this section we use the representation of the ring
$A\cong\F[x]/{\ideal{x^n-1}}$ as in\eqnref{E-A}.
No direct decomposition into fields is needed. Since more than one
automorphism appear simultaneously we do not use the abbreviation $\cR$
for \Azs.
It will be convenient in the following to index
the rows and columns of an $n\times n$-matrix as well as the entries of
$n$-vectors from~$0$ to $n-1$.

We begin with classical circulants. Recall the notation~$\p$
and $\v=\p^{-1}$ from\eqnref{E-p}.

\begin{defi}\label{D-Mg}
For $g=\sum_{i=0}^{n-1}g_i x^i\in A$ define
\[
     M_g:=\begin{bmatrix}g_0&g_1&\ldots&g_{n-2}&g_{n-1}\\
                       g_{n-1}&g_0&\ldots&g_{n-3}&g_{n-2}\\
                        \vdots&\vdots&  &\vdots &\vdots\\
                        g_2&g_3&\ldots&g_0&g_1\\
                        g_1&g_2&\ldots&g_{n-1}&g_0
          \end{bmatrix}
         =\begin{bmatrix} \v(g)\\ \v(xg)\\ \vdots\\ \v(x^{n-2}g)\\
                       \v(x^{n-1}g)\end{bmatrix}
         =\big[g_{(j-i)~\mod n}\big]_{i,j=0,\ldots,n-1}\in\F^{n\times n}.
\]
We call $M_g$ the circulant matrix associated with~$g$.
\end{defi}

The following properties of circulant matrices are either trivial or
well-known in the theory of block codes, see
e.~g.~\cite[p.~501]{MS77}, but also~\cite{Da79} for a general reference on
circulant matrix theory.

\begin{lemma}\label{R-MgMh}
\begin{alphalist}
\item The mapping $A\longrightarrow\F^{n\times n},\ g\longmapsto M_g$
      is $\F$-linear and injective.
\item For $g,h\in A$ we have $M_{gh}=M_gM_h=M_hM_g$.
\item $\rank M_g=\deg\frac{x^n-1}{\gcd(g,x^n-1)}=:k$
      (where the quotient is evaluated in $\F[x]$) and
      every set of $k$ consecutive rows (resp.\ columns) of $M_g$ is linearly
      independent.
\item Let $g=\sum_{i=0}^{n-1} g_i x^i$ and put
      $\widehat{g}:=g(x^{n-1})=g_0+g_{n-1}x+g_{n-2}x^2+\ldots+g_1x^{n-1}$. Then
      \[
              \trans{M_g}=M_{\widehat{g}}.
      \]
      The map $\theta: A\longrightarrow A,\
      g\longmapsto\widehat{g}$ is an involutive $\F$-algebra automorphism
      of~$A$.
\item A matrix $M\in\F^{n\times n}$ is the circulant matrix associated with some
      polynomial $g\in A$ if and only if $SM=MS$ where
      \begin{equation}\label{e-S}
           S:=M_x=\begin{bmatrix}0&1&\cdots  & 0 \\
                 \vdots &\vdots &\ddots&\vdots  \\
                  0 & 0 &\cdots  & 1\\ 1& 0 &\cdots &0 \end{bmatrix}.
      \end{equation}
\item  $ M_g=g(S) \text{ for all }g\in A $.
\item $\v (ab)=\v (a)b(S) $ for all $a,b\in A$.
\item $S $ is the matrix of the
      linear map $a\longmapsto ax$ for $a\in A$ with respect to the basis
      $1,x,\dots,x^{n-1}$ of the $\F$-vector space~$A$ and, as a
      consequence, $M_g=g(S)$ is the matrix of the map $a\longmapsto ag$.
\end{alphalist}
\end{lemma}

The equation in Lemma \ref{R-MgMh}(f) can also be used as an alternative, but less intuitive
definition of circulant matrices. Many of the properties above are
easily proved on the basis of this identity, as there are linearity,
commutativity, and multiplicativity as well as the transposition rule,
where the latter is a direct consequence of the rule $\trans{S}=S^{-1}$.
One also obtains the well known fact that all circulant matrices can be
simultaneously diagonalized over an extension field of~$\F$ that contains a
primitive $n$-th root of unity.

Also for later use we note that the set of all $n\times n$-circulant
matrices over $\F$ is just $\F[S]$ and thus is a commutative subring of
$\F^{n\times n}$ which is isomorphic to~$A$.

The main additional ingredient for our generalized $\sigma$-circulants will be the
following.

\begin{defi}\label{D-Msigma}
For $\sigma\in\AutF(A)$ and $g\in A$ we define
\[
   \Psigma :=\begin{bmatrix} \v(1)\\ \v(\sigma(x))\\ \vdots\\
                  \v(\sigma(x^{n-2}))\\ \v(\sigma(x^{n-1}))
                \end{bmatrix}.
\]
\end{defi}

One should observe that~$\Psigma$ is the matrix with respect to the basis
$1,x,\ldots,x^{n-1}$ associated with the $\F $-linear map which is induced
by the automorphism $\sigma$, i.~e. we have
\begin{equation}\label{e-vPsigma}
  v\Psigma=\v\big(\sigma(\p(v))\big)\text{ for all }v\in\F^n.
\end{equation}

We will need the following properties.

\begin{lemma}\label{L-Msigma}
Let $\sigma,\ \tau\in\AutF(A)$ and $g,\, h\in A$. Then
\begin{arabiclist}
\item $P_{\mbox{\tiny\rm id}}=I_n$ and
      $P_{\sigma\tau}=P_{\tau}P_{\sigma}$.
      Furthermore $\Psigma \in Gl_n(\F)$ and $\Psigma ^{-1}=P_{\sigma^{-1}}$.
\item $P_{\sigma}^{-1}M_g\Psigma =M_{\sigma(g)}$.
\item For $v\in\F^n$ one has $\p (v\Psigma M_g)=\sigma(\p(v)) g$.
\end{arabiclist}
\end{lemma}

\begin{proof}
(1) is a direct consequence of the fact that
$\Psigma,\, P_{\tau}$
are just the matrices which are associated with $\sigma $ and $\tau $
when considered as $\F $-linear maps. The most important
property (2) can be obtained as follows.
For $v\in\F^n $ let $\p (v)=:f$ and suppose $\sigma (x)=a $.
Using Lemma \ref{R-MgMh}(f),~(g) and~(h) as well as\eqnref{e-vPsigma}
we compute
\[
  vS\Psigma = \v (\sigma (fx))=\v (\sigma (f)\sigma (x))=\v
  (\sigma (f))a(S)=v\Psigma M_{\sigma(x)}.
\]
Part~(3) is a direct consequence of the fact that $\Psigma $ and
$M_g $ are matrix representations of the $\F $-linear maps
$\sigma $ and multiplication with $g$.
\end{proof}

Notice that~(2) of the lemma above shows that the automorphisms on~$A$
appear as inner automorphisms
$M_g\mapsto M_{\sigma(g)}=P_{\sigma}^{-1}M_g\Psigma$ on $\F[S]$, where
$\F[S]\cong A$ as noted above.
This observation leads to

\begin{lemma}\label{L-inneraut}
Let $\sigma\in\AutF(A)$ and assume that
\[
   Q^{-1}M_xQ=\Psigma^{-1}M_x\Psigma
\]
for some invertible matrix $Q\in\F^{n\times n}$. Then
$Q=\Psigma M_u$ for some unit $u\in A$.
\end{lemma}

\begin{proof}
By Lemma~\ref{R-MgMh}(e), the identity
$M_x\big(Q\Psigma^{-1}\big)=\big(Q\Psigma^{-1}\big)M_x$
is possible only if $QP_{\sigma}^{-1}$ is a circulant.
Hence $QP_{\sigma}^{-1}=M_{u'}$ for some $u'\in A$ which, by invertibility
of~$Q$, has to be a unit in~$A$.
Using Lemma~\ref{L-Msigma}(2) we obtain
$Q=\Psigma M_{\sigma(u')}$ and $u:=\sigma (u')$ is a unit in
$A$, too.
\end{proof}

Now we can define polynomial circulant matrices.

\begin{defi}\label{D-sigma.circ}
Let $\sigma\in\AutF(A)$. For
$g=\sum_{\nu\geq0}z^{\nu} g_{\nu}\in \Azs$ we define
\[
  \cMsigma{g}:=
  \sum_{\nu\geq0}z^{\nu}\Psigma^{\nu}M_{g_{\nu}}
  \in\F[z]^{n\times n}.
\]
We call $\cMsigma{g}$ the \scirc\ {}(matrix) for~$g$.
\end{defi}

Let us first present an

\begin{exa}\label{E-cMsigma}
Consider again the situation of Example~\ref{E-CCC1}(1) where $\sigma$
is the automorphism given by $\sigma(x)=\alpha^2 x$ and
\[
   g:=(1+\alpha x+\alpha^2 x^2)+z(1+x+x^2)+z^2(1+\alpha^2 x+\alpha x^2).
\]
Then
\[
  P_{\sigma^0}=I_3,\
  \Psigma=\begin{bmatrix}1& &\\ &\alpha^2& \\ & & \alpha\end{bmatrix},\
  P_{\sigma^2}=\begin{bmatrix}1& &\\ &\alpha& \\ & & \alpha^2\end{bmatrix}
\]
and thus
\begin{align*}
  \cMsigma{g}&=\begin{bmatrix}1&\alpha&\alpha^2\\
                 \alpha^2&1&\alpha\\\alpha&\alpha^2&1\end{bmatrix}
                +z\Psigma
                \begin{bmatrix}1&1&1\\1&1&1\\1&1&1\end{bmatrix}
                +z^2P_{\sigma^2}
                \begin{bmatrix}1&\alpha^2&\alpha\\\alpha&1&\alpha^2\\
                  \alpha^2&\alpha&1\end{bmatrix}\\
             &=\begin{bmatrix}1+z+z^2&\alpha+z+\alpha^2z^2&\alpha^2+z+\alpha z^2\\
                \alpha^2+\alpha^2 z+\alpha^2 z^2&1+\alpha^2 z+\alpha z^2&
                   \alpha+\alpha^2z+z^2\\
                \alpha+\alpha z+\alpha z^2&\alpha^2+\alpha z+z^2&
                1+\alpha z+\alpha^2 z^2\end{bmatrix}.
\end{align*}
From Example~\ref{E-CCC1}(1) we conclude that $\cMsigma{g}$ is basic of
rank~$1$. Furthermore, it generates the $1$-dimensional code
$\cC:=\im\cMsigma{g}=
   \im[1+z+z^2,\,\alpha+z+\alpha^2z^2,\,\alpha^2+z+\alpha z^2]
   \subseteq\F^3$.
As noted in~\ref{E-CCC1}(1) the free distance is~$9$.
\mbox{}\hfill$\Box$
\end{exa}


Notice that $\cMsigma{g}=M_g$ whenever $g\in A$.
Hence $\sigma$-circulants form a generalization of
classical circulant matrices. The latter provide a direct link
between ideal generators and generator matrices for cyclic block
codes.
The next proposition shows that a similar link exists
for $\sigma$-circulants and $\sigma$-CCC's.
This will be exploited extensively in Section~\ref{S-gencontr}
where the correspondences  between left principal ideals,
$\sigma$-circulants and $\sigma$-CCC's will be investigated in detail.

\begin{prop}\label{P-cMsigma}
In the situation of Definition \ref{D-sigma.circ} one has
\begin{alphalist}
\item
\[
 \cMsigma{g}
 =\begin{bmatrix}\v(g)\\ \v(xg)\\ \vdots\\ \v(x^{n-1}g)\end{bmatrix}.
\]
\item $\p\big(u\cMsigma{g}\big)=\p(u)g \text{ for all }u\in\F[z]^n$.
In particular,
$\p\big((1,0,\ldots,0)\cMsigma{g}\big)=g$.
\end{alphalist}
\end{prop}
Note that the foregoing rules are equally valid for
classical circulants.

{\sc Proof:}
(a) Let $g=\sum_{\nu\ge 0}z^{\nu}g_{\nu} $.
The $i $-th canonical basis vector in
$\F[z]^n $ is $e_i:=\v (x^i)$. It is sufficient to show that
$\p (e_i\cMsigma{g})=x^ig$ for $1\le i\le n$ .
For this one computes
\[
   \p (e_i\cMsigma{g})=
   \p (\sum_{\nu\ge 0}z^{\nu}e_i \Psigma^{\nu}M_{g_{\nu}})=
   \sum_{\nu\ge 0}z^{\nu}\p
   (e_iP_{\sigma^{\nu}}M_{g_{\nu}})=
   \sum_{\nu\ge 0}z^{\nu}\sigma^{\nu}(x^i)g_{\nu}=
   x^ig,
\]
where we used the $\F[z]$-linearity of $\p$
and Lemma \ref{L-Msigma}(3) in the second and third
equation, respectively\\
(b)  Using~(a) and $\F[z]$-linearity of $\v$ and $\p$ we obtain
for $u=(u_0,\ldots,u_{n-1})\in\F[z]^n$ the identities
$$
  \p\big(u\cMsigma{g}\big)=\p\big(\v\big({\tx\sum_{i=0}^{n-1}}u_ix^ig\big)\big)=\p(u)g.
  \eqno\Box
$$


The following generalizes Lemma~\ref{R-MgMh}(a) and~(b) to $\sigma $-circulants.

\begin{theo}\label{T-sMgsMh}
Let $\sigma\in\AutF(A)$ and $g,\,h\in\Azs$.
\begin{alphalist}
\item The mapping ${\mathcal M}^{\sigma}: \Azs\longrightarrow\F[z]^{n\times n}$
      is $\F$-linear and injective.
\item $\cMsigma{g}\cMsigma{h}=\cMsigma{gh}$.
\end{alphalist}
As a consequence, the Piret algebra $\Azs$ and the ring $\cM^{\sigma}\big(\Azs\big)$ of all
\scirc s are isomorphic as $\F$-algebras.
\end{theo}

{\sc Proof:} (a) is a consequence of the definition of
$\cM^{\sigma} $, the invertibility of
$\Psigma $ and the injectivity of the mapping $g\mapsto M_g$.\\
(b) Let $g=\sum_{\nu\geq0}z^{\nu}g_{\nu}$ and
$h=\sum_{\mu\geq0}z^{\mu}h_{\mu}$. Then
$gh=\sum_{l\geq0}z^l\sum_{\nu+\mu=l}\sigma^{\mu}(g_{\nu})h_{\mu}$ and
using Lemma~\ref{L-Msigma}(2) and Lemma \ref{R-MgMh}(b) we get
$$
  \cMsigma{g}\cMsigma{h}
  =\sum_{l\geq0}z^l\sum_{\nu+\mu=l}P_{\sigma}^{\nu}M_{g_{\nu}}
  P_{\sigma}^{\mu}M_{h_{\mu}}
  =\sum_{l\geq0}z^l\sum_{\nu+\mu=l}
  P_{\sigma}^{\nu+\mu}M_{\sigma^{\mu}(g_{\nu})h_{\mu}}
  =\cMsigma{gh}.
  \eqno\Box
$$

Part~(b) above has the interesting consequence, that each left inverse of
a polynomial~$f$ in $\Azs$ is also a right inverse of~$f$, since this is
the case for the ring $\F[z]^{n\times n}$.

One should observe that the isomorphism $g\mapsto\cMsigma{g}$ induces a
left $\F[z]$-module structure on the set $\cM^{\sigma}\big(\Azs\big)$
which is different from the canonical left $\F[z]$-module structure induced by
$\F[z]^{n\times n}$.
Furthermore, by Definition~\ref{D-sigma.circ} and Lemma~\ref{R-MgMh}(f) we see that
\[
   \cM^{\sigma}(g)=g(z\Psigma,S) \text{ for }
   g(z,x)=\sum_{\nu\geq0}z^{\nu}\sum_{i=0}^{n-1}g_{i\nu}x^i\in\Azs.
\]
Therefore the isomorphism between
$\Azs$ and
$\cM^{\sigma}\big(\Azs\big)$ can also be understood as an
evaluation homomorphism whose image is just $\F[z\Psigma,S]$, a
subring of $\F[z]^{n\times n}$.


Next we turn to transposes of \scirc s.
They will occur later on in the context of control polynomials of
$\sigma$-CCC's.
It turns out that these transposes are in general not \scirc, but rather
$\widehat{\sigma}$-circulant, where
$\widehat{\sigma}\in\AutF(A)$ is such that
$P_{\widehat{\sigma}}=\trans{\Psigma}$. Let us begin with an example.

\begin{exa}\label{E-sigma.tau}
Consider the case $\F=\F_4=\{0,1,\alpha,\alpha^2\}$ and $n=5$. Let
$\sigma\in\AutF(A)$ be given by $\sigma(x)=x^2$.
Then it is easy to see that $\trans{\Psigma }=P_{\widehat{\sigma}}$ where
$\widehat{\sigma}\in\AutF(A)$ is given by $\widehat{\sigma}(x)=x^3$.
Consider now the polynomial
\[
  g:=1+\alpha^2x+\alpha^2x^2+x^3+z(1+x+\alpha^2x^2+\alpha^2x^4)\in\Azs
\]
with associated \scirc\
\[
 \cMsigma{g}=\begin{bmatrix}1+z& \alpha^2+z& \alpha^2+z\alpha^2& 1& z\alpha^2\\
                      0&1+z\alpha^2&\alpha^2+z& \alpha^2+z& 1+z\alpha^2\\
                      1+z& z\alpha^2& 1& \alpha^2+z\alpha^2& \alpha^2+z\\
                      \alpha^2+z\alpha^2& 1+z& z& 1+z\alpha^2& \alpha^2\\
                      \alpha^2+z\alpha^2& \alpha^2& 1+z\alpha^2& z& 1+z\end{bmatrix}.
\]
It is clear that if $\trans{\cMsigma{g}}$ is a circulant matrix at
all, then it is defined by the polynomial given in the first column
of $\cMsigma{g}$. Thus let
\[
  \widehat{g}:=(1+x^2+\alpha^2x^3+\alpha^2x^4)+z(1+x^2+\alpha^2x^3+\alpha^2x^4)\in\Azs.
\]
Then one verifies that
$\cMsigma{\widehat{g}}\not=\trans{\cMsigma{g}}$ but rather
$\cM^{\widehat{\sigma}}(\widehat{g})=\trans{\cMsigma{g}}$.
\\
We will come back to this example in the next sections where we
translate this result into codes and their duals.
\mbox{}\hfill$\Box$
\end{exa}

In order to establish an identity of the type
$\trans{\cM^{\sigma}}(g)=\cM^{\widehat{\sigma}}(\widehat{g})$ for any
automorphism $\sigma $,
we need the existence of an automorphism
$\widehat{\sigma}\in\AutF(A)$ such that
$\trans{\Psigma}=P_{\widehat{\sigma}}$.
In fact, this already implies  the desired identity for the
$\sigma$-circulants since for any $g=\sum_{\nu\geq0}z^{\nu}g_{\nu}$ we
obtain from Definition~\ref{D-sigma.circ}, Lemma~\ref{R-MgMh}(d) and
Lemma~\ref{L-Msigma}(2)
\begin{equation}\label{E-tPsigma}
  \trans{\cM^{\sigma}}(g)
      =\sum_{\nu\geq0}z^{\nu}M_{\widehat{g_{\nu}}}\trans{\Psigma^{\nu}}
     =\sum_{\nu\geq0}z^{\nu}M_{\widehat{g_{\nu}}}P_{\widehat{\sigma}}^{\nu}
      =\sum_{\nu\geq0}z^{\nu}P_{\widehat{\sigma}}^{\nu}M_{\widehat{\sigma}^{\nu}(\widehat{g_{\nu}})}
      =\cM^{\widehat{\sigma}}(\wh{g}{\sigma}),
\end{equation}
where
$\wh{g}{\sigma}=\sum_{\nu\geq0}z^{\nu}\widehat{\sigma}^{\nu}(\widehat{g_{\nu}})$.
In order to show the existence of $\widehat{\sigma}$, we will make
use of the involution $\theta$ given in Lemma~\ref{R-MgMh}(d).
Notice that by Lemma~\ref{L-Msigma}(2) and Lemma~\ref{R-MgMh}(d) we
have
\[
  \trans{P_{\sigma}^{-1}}M_x\trans{\Psigma}
  =\trans{(\Psigma M_{\widehat{x}}P_{\sigma^{-1}})}=\trans{M_{\sigma^{-1}(\widehat{x})}}
  =M_{\widehat{\sigma^{-1}(\widehat{x})}}.
\]
Taking into account once more Lemma~\ref{L-Msigma}(2), this indicates how the desired automorphism
$\widehat{\sigma}$ has to look like.

\begin{theo}\label{T-tPsigma}
Let $\sigma\in\AutF(A)$. Define
$\widehat{\sigma}:=\theta\,\circ\,\sigma^{-1}\circ\theta\in\AutF(A)$, thus
$\widehat{\sigma}(a)=\widehat{\sigma^{-1}(\widehat{a})}$ for all $a\in A$.
Then $\widehat{\widehat{\sigma}}=\sigma$ and
$\sigma\longmapsto\widehat{\sigma}$ defines an anti-automorphism on the
group $\AutF(A)$. Furthermore,
\[
    \trans{\Psigma}=P_{\widehat{\sigma}}.
\]
\end{theo}

\begin{proof}
Only the identity $\trans{\Psigma}=P_{\widehat{\sigma}}$ needs proof.
Applying several times Lemma~\ref{L-Msigma}(2)
and Lemma~\ref{R-MgMh}(d)
we obtain
\[
  P_{\widehat{\sigma}}^{-1}M_xP_{\widehat{\sigma}}=
  M_{\widehat{\sigma}(x)}=\trans{M_{\widehat{\widehat{\sigma}(x)}}}
  =\trans{\big(\Psigma M_{\sigma(\widehat{\widehat{\sigma}(x)})}P_{\sigma}^{-1}\big)}
  =\trans{P_{\sigma}^{-1}}M_{\widehat{\sigma(\widehat{\widehat{\sigma}(x)})}}
  \trans{\Psigma}
  =\trans{P_{\sigma}^{-1}}M_x\trans{\Psigma},
\]
where the last equality follows from
$\widehat{\sigma(\widehat{\widehat{\sigma}(x)})}=x$.
Lemma~\ref{L-inneraut} now yields
\begin{equation}\label{E-tPsigma2}
   \trans{\Psigma}=P_{\widehat{\sigma}}M_u\text{ for some unit }u\in A.
\end{equation}
We will show now that the matrix $\Psigma$ not only has zeros in the
first row except for the very first entry (which is obvious by definition),
but also in the first column.
Then Equation\eqnref{E-tPsigma2} implies that the first row of~$M_u$ is of
the form $(1,0,\ldots,0)$ and, being a circulant, $M_u=I_n$. This proves
the theorem.
\\
In order to establish the zero entries in the first column of $\Psigma$
let $\sigma(x)=a=\sum_{l=0}^{n-1}a_lx^l$.
For the rest of the proof it will be convenient to use the notation
$[f]_i$ for the coefficient of~$x^i$ in the polynomial $f\in A$.
Then, according to Definition~\ref{D-Msigma}
we have to show $[\sigma(x^i)]_0=[a^i]_0=0$ for all $i>0$.
Since $\sigma$ is an automorphism, the powers $1,\,a,\ldots,a^{n-1}$
are linearly independent over $\F$ and $a^n=1$. Using linearity and
multiplicativity of the circulants, this implies that the
characteristic polynomial of $M_a$ is given by $X^n-1$ and we can
conclude\quad $0=\mbox{\rm trace}(M_{a})=n [a]_0$. But then also
$[a]_0=0$ since $\gcd(n,\mbox{\rm char}(\F))=1$. As for $[a^i]_0$, we
wish to argue along the same lines. In order to determine the
characteristic polynomial of $M_{a^i}=(M_a)^i$, let
$X^n-1=\prod_{l=0}^{n-1}(X-\omega^l)$ for some
primitive $n$-th root of unity~$\omega$ in some extension field
$\overline{\F}$ of~$\F$. Furthermore assume $\gcd(i,n)=d$ and
$n=d\tilde{n}$. Then
$\omega^i$ is a primitive $\tilde{n}$-th root of unity and, since
$M_a$ is diagonalizable over~$\overline{\F}$, the characteristic
polynomial of $M_{a^i}$ is given by
\[
  \prod_{l=0}^{n-1}(X-\omega^{il})
  =\Big(\prod_{l=0}^{\tilde{n}-1}(X-\omega^{il})\Big)^d
  =(X^{\tilde{n}}-1)^d.
\]
As above we conclude $n[a^i]_0=\mbox{\rm trace}(M_{a^i})=0$ for all
$i>0$ (in which case $\tilde{n}>1$) and, again, $[a^i]_0=0$.
\end{proof}

In Example~\ref{E-sigma.tau} above we had the specific situation that
$\widehat{\sigma}=\sigma^{-1}$. This need not be the case in
general, see the remark below.
However, the automorphisms of $A=\F_4[x]/\ideal{x^5-1}$ as listed in
Example~\ref{E-autos}(b) all satisfy either
$\widehat{\sigma}=\sigma$ (in which case $\Psigma$ is symmetric) or
$\widehat{\sigma}=\sigma^{-1}$ (which only occurs for
$\sigma(x)=x^2$ and $\sigma(x)=x^3$).
This too, is not true in general.

\begin{rem}\label{R-sigmahat}
For the special class of automorphisms~$\sigma$ satisfying
\begin{equation}\label{e-axr}
     \sigma(x)=\gamma x^r,
\end{equation}
where $\gamma\in\F,\,r\in\{1,\ldots,n-1\}$, the associated
$\widehat{\sigma}$ can be found easily.
First notice that\eqnref{e-axr}
induces an automorphism $\sigma\in\AutF(A)$ if and only if
$\gamma^n=1$ and $\gcd(r,n)=1$.
If these conditions are satisfied, then
$\sigma^{-1}$ and $\widehat{\sigma}$ are given by
the equations
\begin{equation}\label{e-sigmahat}
   \sigma^{-1}(x)=\gamma^{-l}x^l\text{ and }
   \widehat{\sigma}(x)=\gamma^lx^l\text{ where }lr\equiv 1~\mod n.
\end{equation}
This can be verified remembering the definition of
$\widehat{\sigma}$.
The conditions in\eqnref{e-sigmahat} lead to plenty of examples where the automorphisms
$\sigma,\,\sigma^{-1},\,\widehat{\sigma}$ are all different,
e.~g. for $A=\F_4[x]/\ideal{x^7-1}$ and~$\sigma$ given by
$\sigma(x)=\alpha x^4$.
\\
We also wish to note that in~\cite{Pi76} only automorphisms
as in\eqnref{e-axr} with $\gamma=1$ were considered.
In this case one always has $\widehat{\sigma}=\sigma^{-1}$.
\end{rem}

Now we can describe the transposes of \scirc s. In part~(a) below we
obtain a direct generalization of Lemma~\ref{R-MgMh}(d). The
anti-isomorphism in part~(b) will be crucial in the next section when
relating a control polynomial of a $\sigma$-cyclic code~$\cC$ to a
generator polynomial of the $\widehat{\sigma}$-cyclic dual code $\cC^{\perp}$.
\begin{theo}\label{T-Msigmatransp}
Let $\sigma\in\AutF(A)$ and $\widehat{\sigma}$ be defined as in Theorem~\ref{T-tPsigma}.
For any polynomial
$g=\sum_{\nu\geq0}z^{\nu}g_{\nu}\in\Azs$ define
\begin{equation}\label{E-tM2}
  \wh{g}{\sigma}:=\sum_{\nu \ge 0}\widehat{g_{\nu}}z^{\nu}=
  \sum_{\nu \ge 0}z^{\nu}\widehat{\sigma}^{\nu}(\widehat{g_{\nu}}) \
  \in A[z;\widehat{\sigma}].
\end{equation}
Then
\begin{alphalist}
\item $\trans{\cMsigma{g}}=\cM^{\widehat{\sigma}}(\wh{g}{\sigma})$.
\item The map
      $\wh{\ \ }{\sigma} :\ \Azs\longrightarrow A[z;\widehat{\sigma}],\ g\longmapsto\wh{g}{\sigma}$ is
      an anti-isomorphism of the $\F$-algebras $\Azs$ and $A[z;\widehat{\sigma}]$,
      that is, $\wh{\ \ }{\sigma}$ is $\F$-linear and
      satisfies $\wh{gh}{\sigma}=\wh{h}{\sigma}\wh{g}{\sigma}$ for all
      $g,\,h\in\Azs$.
      The inverse map is given by $\wh{\ \ }{\widehat{\sigma}} :
      \ A[z;\widehat{\sigma}] \longrightarrow A[z;\sigma],\ g\longmapsto\wh{g}{\widehat{\sigma}} $.
\end{alphalist}
\end{theo}

\begin{proof}
(a) has been shown in\eqnref{E-tPsigma}.
\\
(b) $\F$-linearity and injectivity are obvious by\eqnref{E-tM2}.
Anti-multiplicativity is a consequence of
\[
  \cM^{\widehat{\sigma}}(\wh{gh\,}{\sigma})=\trans{\cMsigma{gh}}
  =\trans{\cMsigma{h}}\trans{\cMsigma{g}}
  =\cM^{\widehat{\sigma}}(\wh{h}{\sigma})\cM^{\widehat{\sigma}}(\wh{g}{\sigma})
  =\cM^{\widehat{\sigma}}(\wh{h}{\sigma}\wh{g}{\sigma})
\]
along with injectivity of the map $\cM^{\widehat{\sigma}}$.
Finally, the equation
$\cMsigma{g}=\trans{\cM^{\widehat{\sigma}}(\wh{g}{\sigma})}=
 \cMsigma{\wh{(\wh{g}{\sigma})}{\widehat{\!\!\sigma}}}$ shows that
$\wh{(\wh{g}{\sigma})}{\widehat{\!\!\sigma}}=g$ for all $g\in\Azs$, which completes the
proof.
\end{proof}

As a simple consequence of Theorem~\ref{T-Msigmatransp} we obtain
that each polynomial vector appears
not only as a row but also as a column in some
\scirc. Furthermore, as we will show next, the algebra of \scirc s is
saturated in the sense that if a multiple of a circulant within the
ring~$\F[z]^{n\times n}$ is a circulant again, then it is even a multiple
within the algebra $\cMsigma{\Azs}$.
Also these results will be of use in the next section
for generator and control matrices of $\sigma$-CCC's.


\begin{cor}\label{C-circs}
Let $\sigma\in\AutF(A)$. Then one has the implications
\begin{arabiclist}
\item For each $v\in\F[z]^n$ and
$g_1=\p(v)\in A[z;\sigma],\,g_2=\p(v)\in A[z;\widehat{\sigma}]$
and for each $f\in\Azs$ one has
\[
          v\cMsigma{f}=0\, \Longleftrightarrow\,
          \cMsigma{g_1}\cMsigma{f}=0\ \text{ and }\
       \cMsigma{f}\trans{v}=0\, \Longleftrightarrow\,
          \cMsigma{f}\cMsigma{\wh{g_2}{\widehat{\sigma}}}=0.
\]
\item For all $f,g\in\Azs$ one has the two implications
\begin{align*}
   &\exists\;Q\in\F[z]^{n\times n}:\cMsigma{f}=Q\cMsigma{g}
   \Longrightarrow\exists\;h\in\Azs:\cMsigma{f}=\cMsigma{h}\cMsigma{g},
   \\[.5ex]
  &\exists\;Q\in\F[z]^{n\times n}:\cMsigma{f}=\cMsigma{g}Q
   \Longrightarrow\exists\;h\in\Azs:\cMsigma{f}=\cMsigma{g}\cMsigma{h}.
\end{align*}
\end{arabiclist}
\end{cor}

One should observe that part~(2) above, applied to constant polynomials
$f,\, g\in A$ leads to the analogous statements for classical circulants.

\begin{proof}
(1)\, The first equivalence can
be obtained as follows with the help of Proposition
\ref{P-cMsigma}(b):
\[
v\cMsigma{f}=0\Longleftrightarrow
g_1f=0\Longleftrightarrow
\cMsigma{g_1}\cMsigma{f}=0.
\]
The second equivalence follows from the first one by transposition
and Theorem \ref{T-Msigmatransp}.  \\
(2)\, Let $v$ be the first row of $\cMsigma{f} $ and $w$ be
the first row of $Q$ and $h=\p (w) $.
Then $\v (f)=v $ and $v=w\cMsigma{g}$. By Proposition
\ref{P-cMsigma}(b) we obtain $f=hg $ which gives us
$\cMsigma{f}=\cMsigma{h}\cMsigma{g} $. The second statement
follows as in (a) by transposition and Theorem \ref{T-Msigmatransp}.
\end{proof}

So far we have not discussed the rank of \scirc s. As opposed to classical
circulants (see Lemma~\ref{R-MgMh}(c)) there is no general simple rule telling the
rank of $\cMsigma{g}$ based on the polynomial~$g$. Fortunately, if~$g$ is
a reduced polynomial, a generalization of the classical result exists. This
will be treated in Theorem~\ref{P-rank}.

\section{Description of $\sigma $-cyclic codes and their
duals}\label{S-gencontr}
\setcounter{equation}{0}

Now we are in a position to return to $\sigma$-CCC's in the sense
of Definition~\ref{D-Roos} or Observation~\ref{Piretalgebra}(b).
In this section we introduce generator and control polynomials as
well as (square circulant) generating and control matrices for
$\sigma$-CCC's. We show that they behave just like those for block
codes. Below we first summarize the relation between cyclic block
codes and classical circulant matrices, as this shows exactly what
we are after for convolutional codes. As a reference on cyclic
block codes any (introductory) book on coding theory suffices, for
instance~\cite{MS77} or~\cite{Be98}.

Let $\cC \in\F^n $ be a cyclic block code, then  --- in
polynomial representation --- we obtain a principal left ideal
$\cJ=\p (\cC)=\lideal{g}$ for some $g\in A$.
Once given a generator polynomial $g$, then the classical
circulant $M_g$ is a generating matrix for $\cC$ in the sense of
Proposition~\ref{P-basicLA} and one has
\begin{equation}\label{e-cycbl1}
     \cC:=\v(\cJ)=\im M_g=\ker M_h,
\end{equation}
where $h\in A$ generates the annihilator ideal of $\cJ$ in~$A$ and~$M_h$ is its circulant.

Usually, $M_g$ is not an encoder for~$\cC$, which must have full rank.
Such an encoder is obtained by extracting the first~$k$ rows of~$M_g$,
where $k=\dim_{\F}\cC=\rank M_g$.
Of course, the generator polynomial~$g$ is not unique and
can be modified by multiplying with units~$u$ from~$A$.
For the circulants this amounts to multiplying by~$M_u$
from either side since classical circulants commute.
There are two natural ways of choosing a specific~$g$ by
imposing one of the conditions
\begin{equation}\label{E-norm1}
g\,  \big| \,x^n-1\ \text{ in } \F[x] \text{  and the
leading coefficient is 1}
\end{equation}
or
\begin{equation}\label{E-norm2}
g\, \text{ is idempotent}.
\end{equation}
The first condition is more widely used and the name
'generator polynomial' usually refers to this choice.
If both,~$g$ and~$h$ of\eqnref{e-cycbl1} satisfy\eqnref{E-norm1}
then $gh=x^n-1$ and~$h$ is just the complementing factor for~$g$ and this is
what usually is meant when calling~$h$ a 'control polynomial'.

In the situation of\eqnref{e-cycbl1} the dual code
$\cC^{\perp}:=\{w\in\F^n\mid w\,\trans{v}=0\text{ for all }v\in\cC\}$
is given by
\begin{equation}\label{e-cycbl2}
   \cC^{\perp}=\im\trans{M_h}=\ker\trans{M_g}
   =\ker M_{\widehat{g}}=\im M_{\widehat{h}},
\end{equation}
where $\widehat{g}$ and $\widehat{h}$
are defined as in Lemma~\ref{R-MgMh}(d). Normalizing
according to\eqnref{E-norm1} leads to the polynomials
$h(0)^{-1}x^k\widehat{h}$ (resp.\ $g(0)^{-1}x^{n-k}\widehat{g}$),
the generator (resp.\ control) polynomial of $\cC^{\perp}$,
see e.g.~\cite[p.~196]{MS77}.
Here $h(0)$ and $g(0)$ denote the constant terms of~$h$ and~$g$.

In this section we will show, that with the help of
$\sigma $-circulants and the generator polynomials from
Section~\ref{S-idealgen} and~\ref{gencomp} the complete
scenario generalizes nicely to $\sigma$-CCC's.
In addition, the basic notions of convolutional coding theory, like
non-catastrophicity, minimality, and complexity, can be incorporated
successfully.

Throughout this section let $\sigma\in\AutF(A)$ be a fixed
automorphism and, as before, let $\cR:=\Azs$.

Recall from Observation~\ref{Piretalgebra}(b) that a submodule
$\cC\subseteq\F[z]^n$ is called $\sigma$-cyclic if~$\p(\cC)$ is a left ideal in~$\cR$.
Using the calculus of \scirc s, this can also be expressed in terms of vector
polynomials. One simply has to translate multiplication by~$x$
in~$\cR$ via the isomorphism~$\v$ into a suitable
mapping~$\mathfrak{m}$ on~$\F[z]^n$. Observe that, due to
noncommutativity of~$\cR$, this mapping is $\F$-linear but not
$\F[z]$-linear.

\begin{obs}\label{P-cMsigma2}
A submodule $\cC\subseteq\F[z]^n$ is $\sigma$-cyclic iff
$\mathfrak{m}(\cC)\subseteq\cC$,
where
\[
   \mathfrak{m}:\F[z]^n\longrightarrow\F[z]^n,\quad
   \sum_{\nu\geq0}z^{\nu}v_{\nu}\longmapsto
   \sum_{\nu\geq0}z^{\nu}v_{\nu}M_{\sigma^{\nu}(x)}
   =\sum_{\nu\geq0}z^{\nu}v_{\nu}\Psigma^{-\nu}S\Psigma^{\nu}
\]
and $S=M_x$, as in\eqnref{e-S}. This follows from the fact that
$\mathfrak{m}(v)=\v(x\p(v))$ for all $v\in\cC$, which itself is
equivalent to $\p(\mathfrak{m}(v))=x\p(v)$ and this is a direct
consequence of Proposition~\ref{P-cMsigma}(b) and
Lemma~\ref{L-Msigma}(2).
\end{obs}

Observe that for $\sigma=\mbox{\rm id}$ one has
$M_{\sigma^{\nu}(x)}=S$ and $\Psigma=I_n$ so that in this case~$\mathfrak{m}$
describes the classical cyclic shift.
Furthermore, if $\sigma(x)=x^m$ for some~$m$ that is coprime with~$n$,
then $M_{\sigma^{\nu}(x)}=M_{x^{(m^{\nu})}}=S^{(m^{\nu})}$ and one
obtains\eqnref{Piret-def1}.

By Theorem~\ref{existence} each delay-free $\sigma$-cyclic submodule is a
principal left ideal when considered in $\Azs$. Using the
correspondence of \scirc s and principal left ideals as described in
Proposition~\ref{P-cMsigma}(b) this immediately leads to a circulant
generating matrix. Precisely, one has
\begin{equation}\label{e-ideal.circ}
       \v\big(\lideal{g}\big)=\im\cMsigma{g}\text{ for all }g\in\cR.
\end{equation}
As a consequence, a delay-free submodule~$\cC\subseteq\F[z]^n$ is
$\sigma$-cyclic if and only if $\cC=\im\cMsigma{g}$ for some $g\in\Azs$, which,
additionally, can be taken as a reduced and normalized polynomial
satisfying\eqnref{E-existence}, see Corollary~\ref{C-redu}.

In order to also get a description of $\sigma$-cyclic codes by
control polynomials and control matrices we need the following.

\begin{defi}\label{D-annih}
Let $F\subseteq\cR$ be any subset. Then
\begin{arabiclist}
\item $F^{^{\circ}}:=\{h\in\cR\mid\forall\; f\in F:  fh=0\}$.
\item $^{^{\circ}}\!F:=\{g\in\cR\mid\forall\; f\in F:  gf=0\}$.
\end{arabiclist}
We call $F^{^{\circ}}$ and $^{^{\circ}}\!F$ the right and left annihilator of the
set~$F$, respectively.
Obviously, $F^{^{\circ}}=\rcirc{\lideal{F}}$ and $^{^{\circ}}\!F=\lcirc{\rideal{F}}$ are
the right and left annihilator of the left and right ideal generated by~$F$,
respectively.
\end{defi}

Using Observation~\ref{O-attrib} and the fact that $\F[z]$ does not contain any zero
divisors of~$\cR$, one verifies straightforwardly the following.

\begin{obs}\label{R-annih}
The annihilators $^{^{\circ}}\!F$ and $F^{^{\circ}}$ are direct
summands of~$\cR$ as left resp.\ right $\F[z]$-modules.
In particular, both ideals are delay-free and by
Theorem~\ref{existence} and Corollary~\ref{c-anti} are principal
left resp.\ right ideals.
\end{obs}

Now we have the following

\begin{lemma}\label{L-annih}
Let $g,\,h\in\cR$. Then
\begin{arabiclist}
\item $\rcirc{\lideal{g}}=\rideal{h}\Longleftrightarrow
       \ker\trans{\cMsigma{g}}=\im\trans{\cMsigma{h}}$.
\item $\lideal{g}=\lcirc{\rideal{h}}\Longleftrightarrow
       \im\cMsigma{g}=\ker\cMsigma{h}$.
\end{arabiclist}
Furthermore, if the identities in equivalence in~(1) (resp.~(2)) are satisfied, then the
matrix $\cMsigma{h}$ (resp.\ $\cMsigma{g}$) is basic.
\end{lemma}

\begin{proof}
(1) can be established as follows.
\begin{align*}
  &\rcirc{\lideal{g}}=\rideal{h}\\
  &\quad\Longleftrightarrow gh=0\text{ and }
      [gf=0\Longrightarrow\exists\,a\in\cR: f=ha]\\
  &\quad\Longleftrightarrow\cMsigma{g}\cMsigma{h}=0\text{ and }
     [\cMsigma{g}\cMsigma{f}=0\Longrightarrow\exists\,a\in\cR:
       \cMsigma{f}=\cMsigma{h}\cMsigma{a}]\\
  &\quad\Longleftrightarrow\trans{\cMsigma{h}}\trans{\cMsigma{g}}=0\;\text{and}\;
     [\trans{\cMsigma{f}}\trans{\cMsigma{g}}=0\Rightarrow\exists\,a\in\cR:
       \trans{\cMsigma{f}}=\trans{\cMsigma{a}}\trans{\cMsigma{h}}].
\end{align*}
The last statement is satisfies if and only if
$\ker\trans{\cMsigma{g}}=\im\trans{\cMsigma{h}}$, which can be seen as follows.
\\
For the if-part only the implication in brackets needs proof.
But this is obtained from Corollary~\ref{C-circs}(2) since
$\trans{\cMsigma{f}}\trans{\cMsigma{g}}=0$ along with the assumption implies
$\im\trans{\cMsigma{f}}\subseteq\im\trans{\cMsigma{h}}$, hence
$\cMsigma{f}=\cMsigma{h}Q$ for some matrix~$Q$ and the corollary applies.
\\
For the only-if-part we have to show that
$\ker\trans{\cMsigma{g}}\subseteq\im\trans{\cMsigma{h}}$.
Thus let $v\trans{\cMsigma{g}}=0$ for some $v\in\F[z]^n$.
By Corollary~\ref{C-circs}(1) we obtain
$\cMsigma{g}\cM^{\sigma}(\wh{f}{\widehat{\sigma}})=0$, where
$f=\p(v)\in A[z;\widehat{\sigma}]$.
Then the assumption implies that
$\cM^{\widehat{\sigma}}(f)=\trans{\cMsigma{\wh{f}{\widehat{\sigma}}}}=
\trans{\cMsigma{a}}\trans{\cMsigma{h}}$ for some $a\in\Azs$ and hence
$v\in\im\trans{\cMsigma{h}}$.
\\
(2) In this case the anti-isomorphism $\wh{\ }{\sigma}$ of
Theorem~\ref{T-Msigmatransp}(b) yields
$\lideal{g}=\lcirc{\rideal{h}}$ in $\Azs$ if and only if
$\rideal{\wh{g}{\sigma}}=\rcirc{\lideal{\wh{h}{\sigma}}}$ in $A[z;\widehat{\sigma}]$.
Thus, use of~(1) and Theorem~\ref{T-Msigmatransp}(a) leads to the desired result.
\\
The additional assertion that the two given matrices are basic follows
either from the equivalence of Proposition~\ref{P-directsummand}(6) and~(3) or
from the direct summand property as stated in Observation~\ref{R-annih} together
with~\ref{P-directsummand}(5).
\end{proof}

The following theorem collects the basic facts on $\sigma$-CCC's.
Recall that transposes of \scirc s are $\widehat{\sigma}$-circulants.
Therefore, the dual code of a $\sigma$-CCC corresponds to a left ideal in
the Piret algebra $A[z;\widehat{\sigma}]$.
For simplicity, we use the notation $\lideal{\ }$ for left ideals in
either Piret algebra; yet, in order to avoid confusion, we will make
the corresponding algebra precise at each point.
Recall also that the isomorphism~$\p$ in\eqnref{E-p} does not depend on the
multiplicative structure of the set~$A[z]$ so that we may use it for
either Piret algebra.

\begin{theo}\label{T-annihMat}
Let $\cC\subseteq\F[z]^n$ be a $\sigma$-cyclic code and let
\[
     \cC^{\perp}:=\big\{w\in\F[z]^n\,\big|\, w\trans{v}=0\text{ for all }v\in\cC\big\}
\]
be its dual code.
Furthermore, let $g,\,h\in A[z;\sigma]$ be such that $\p(\cC)=\lideal{g}$ and
$\rcirc{\lideal{g}}=\rideal{h}$.
Then
\begin{alphalist}
\item $\cMsigma{g}$ and $\cMsigma{h}$ are both basic.
\item $\cC=\im\cMsigma{g}=\ker\cMsigma{h}$.
\item $\p(\cC)=\lideal{g}=\lcirc{\rideal{h}}$ in $\Azs$.
\item $\cC^{\perp}=\ker\cM^{\widehat{\sigma}}(\wh{g}{\sigma})
       =\im\cM^{\widehat{\sigma}}(\wh{h}{\sigma})$.
\item $\p(\cC^{\perp})=\lideal{\wh{h}{\sigma}}$ in
      the Piret algebra $A[z;\widehat{\sigma}]$. Hence the
      dual of a $\sigma$-CCC is a $\widehat{\sigma}$-CCC.
\end{alphalist}
\end{theo}

\begin{proof}
(a) $\cMsigma{g}$ is basic since it generates a code, see
Proposition~\ref{P-directsummand}(3); $\cMsigma{h}$ is basic by
Lemma~\ref{L-annih}.
\\
(b) $\cC=\im\cMsigma{g}\subseteq\ker\cMsigma{h}$ follows from the choice of~$g$
and~$h$, see also\eqnref{e-ideal.circ}.
Furthermore, by Lemma~\ref{L-annih}(1),
$\ker\trans{\cMsigma{g}}=\im\trans{\cMsigma{h}}$ and thus
$\rank\cMsigma{h}=n-\rank\cMsigma{g}$.
But then Proposition~\ref{P-directsummand}(7) yields
$\cC=\ker\cMsigma{h}$ since~$\cC$ is a direct summand.
\\
(c) is a consequence of~(b) along with Lemma~\ref{L-annih}(2).
\\
(d) follows from the obvious fact that
$\cC^{\perp}=\ker\trans{\cMsigma{g}}=\im\trans{\cMsigma{h}}$.
\\
(e) is a consequence of~(d) along with Lemma~\ref{L-annih}(1).
\end{proof}

These results motivate the following definition.

\begin{defi}\label{D-gencontr}
Let $\cC\subseteq\F[z]^n$ be a $\sigma$-cyclic code and $g,\,h\in\cR$
be such that $\p(\cC)=\lideal{g}$ and $\p(\cC)^{\circ}=\rideal{h}$.
Then we call~$g$ a generator polynomial and~$h$ a control polynomial of the
code~$\cC$.
Consequently, the polynomials $\wh{h}{\sigma}$ and $\wh{g}{\sigma}\in A[z;\widehat{\sigma}]$
are a generator and a control polynomial of the dual code~$\cC^{\perp}$, respectively.
\end{defi}

At this point there is no need to normalize generator and control
polynomials.
But there is a way to obtain uniqueness by requiring~$g$ and~$\wh{h}{\sigma}$ to be
left reduced and their $z$-free terms to be normalized according
to\eqnref{E-norm1} or\eqnref{E-norm2}.

Via the anti-isomorphism
$\wh{\ \ }{\sigma} :\ \Azs\longrightarrow A[z;\widehat{\sigma}],\
g\longmapsto\wh{g}{\sigma}$ from Theorem~\ref{T-Msigmatransp},
one observes that the right annihilator $\rideal{h}$ of the code
$\cC=\v(\lideal{g})$ is anti-isomorphic to the dual code
$\cC^{\perp}=\v(\lideal{\wh{h}{\sigma}})$.

The following very detailed example is designed to shed some light on all aspects
of our setting thus far.

\begin{exa}\label{E-sigma.tau2}
Consider again Example~\ref{E-sigma.tau} where $\F=\F_4,\,n=5$, and
$\sigma(x)=x^2$.
The circulant $\cMsigma{g}$ associated with the polynomial
\[
  g:=1+\alpha^2x+\alpha^2x^2+x^3+z(1+x+\alpha^2x^2+\alpha^2x^4)\in\Azs
\]
can be shown to be basic. Since $\rank\cMsigma{g}=2$, it defines a
$2$-dimensional $\sigma$-cyclic code $\cC\subseteq\F[z]^5$.
A control polynomial, i.~e. a right annihilator of the left ideal $\lideal{g}\in\Azs$
can be found as follows.
First we compute a basis $w_1,\,w_2,\,w_3\in\F[z]^5$ of the right kernel
of $\cMsigma{g}$, i.~e. $\cMsigma{g}\trans{w_i}=0$ for $i=1,2,3$.
This can easily be achieved by use of a Smith-form of $\cMsigma{g}$ and
yields the basic matrix
\[
   H:=\big[\trans{w_1},\trans{w_2},\trans{w_3}]=
   \begin{bmatrix}1& 0& 0 \\ 0& 1& 0 \\1+z\alpha^2& z\alpha^2& \alpha+z\alpha+z^2 \\
                  \alpha+z& \alpha^2+z&1+z+z^2\alpha \\
                  \alpha+z\alpha^2&\alpha^2+z\alpha^2&\alpha+z^2\end{bmatrix}.
\]
Then $\im\trans{H}=\ker\cM^{\widehat{\sigma}}(\wh{g}{\sigma})$ and any
control polynomial $h\in\Azs$ of~$\cC$ satisfies
$\im\trans{H}=\im\trans{\cMsigma{h}}=\im\cM^{\widehat{\sigma}}(h')$ where $h':=\wh{h}{\sigma}$.
In $A[z;\widehat{\sigma}]$ this reads as $\lideal{\p(w_1),\,\p(w_2),\,\p(w_3)}=\lideal{h'}$
and thus we need to find a principal generator of
this left ideal in $A[z;\widehat{\sigma}]$.
Hence put
\begin{align*}
  f_1:=&\p(\trans{w_1})=1+x^2+\alpha x^3+\alpha x^4+z(\alpha^2 x^2+x^3+\alpha^2 x^4),\\
  f_2:=&\p(\trans{w_2})=x+\alpha^2x^3+\alpha^2x^4+z(\alpha^2 x^2+x^3+\alpha^2 x^4),\\
  f_3:=&\p(\trans{w_3})=\alpha x^2+x^3+\alpha x^4+z(\alpha x^2+x^3)
               +z^2(x^2+\alpha x^3+x^4).
\end{align*}
By Observation~\ref{R-annih} we know that $\rideal{h}$ is delay-free, thus,
using the anti-automorphism between $\Azs$ and $A[z;\widehat{\sigma}]$, we
get the delay-freeness of the left ideal $\lideal{h'}$ in
$A[z;\widehat{\sigma}]$.
As a consequence, application of Algorithm~\ref{Algo} to the family
$f_1,\,f_2,\,f_3$ produces the desired principal generator~$h'\in A[z;\widehat{\sigma}]$,
even in reduced form.
In order to actually perform these computation we first need to know the
automorphism $\widehat{\sigma}$. Using Remark~\ref{R-sigmahat} we find
$\widehat{\sigma}(x)=x^3$.
Furthermore, the algorithm needs a decomposition of $A\cong\F[x]/{\ideal{x^5-1}}$
into a direct sum of fields and the representation of the
automorphism~$\widehat{\sigma}$ as well as the given data in the according
form.
For this task we may use the list in Example~\ref{E-autos}(b).
Switching to the notation of the second column therein we find
$\widehat{\sigma}[a,b,c]=[a,\Psi^{-1}(c),\Psi(b)^4]$ where~$\Psi$ is as
in\eqnref{e-psi}.
Furthermore, using the Chinese Remainder Theorem, precisely the map~$\varrho$
given in\eqnref{E-rho}, the polynomials $f_1,\,f_2,\,f_3$ turn into
\begin{align*}
 \varrho(f_1)&=h_1:=\ve{2}\alpha^2x+z\ve{1}+z\ve{3}(\alpha x+\alpha),\\
 \varrho(f_2)&=h_2:=\ve{1}+\ve{2}\alpha^2x+z\ve{1}+z\ve{3}(\alpha x+\alpha),\\
 \varrho(f_3)&=h_3:=\ve{1}+\ve{2}(\alpha^2 x+\alpha^2)+z\ve{1}\alpha^2+
                         z\ve{2}x+z\ve{3}(\alpha x+1)\\
               &\qquad\qquad +z^2\ve{1}\alpha+z^2\ve{3}(\alpha^2 x+\alpha^2),
\end{align*}
where $\ve{1}=[1,0,0],\,\ve{2}=[0,1,0],\,\ve{3}=[0,0,1]$.
Now we may run Algorithm~\ref{Algo} on the data $h_1,h_2,h_3$.
Exactly this has been done in Example~\ref{E-comp}.
Therein, a principal generator of the left ideal $\lideal{h_1,h_2,h_3}$ in
$A[z;\widehat{\sigma}]$ was found to be the reduced and normalized polynomial
\begin{equation}\label{e-h2prime}
   h'':=\ve{1}+\ve{2}+z\ve{3}.
\end{equation}
Translating this back we obtain
\begin{equation}\label{e-hprime}
   h':=\varrho^{-1}(h'')=1+\alpha^2 x+\alpha x^2+\alpha x^3+\alpha^2x^4+
                 z(\alpha^2x+\alpha x^2+\alpha x^3+\alpha^2x^4)
\end{equation}
and finally the control polynomial
\begin{equation}\label{e-hpoly}
  h=\wh{h'}{\widehat{\sigma}}=
  1+\alpha^2x+\alpha x^2+\alpha x^3+\alpha^2x^4
  +z(\alpha x+\alpha^2x^2+\alpha^2x^3+\alpha x^4)
\end{equation}
of the code~$\cC$ as well as the associated circulant control matrix
\begin{equation}\label{e-Msigmah}
  \cMsigma{h}=\begin{bmatrix}1&
        \alpha^2+z\alpha & \alpha +z\alpha^2& \alpha +z\alpha^2& \alpha^2+z\alpha \\
        \alpha^2+z\alpha^2& 1+z\alpha & \alpha^2& \alpha +z\alpha & \alpha +z\alpha^2\\
        \alpha +z\alpha & \alpha^2+z\alpha^2& 1+z\alpha^2& \alpha^2+z\alpha & \alpha \\
        \alpha +z\alpha & \alpha & \alpha^2+z\alpha & 1+z\alpha^2& \alpha^2+z\alpha^2\\
        \alpha^2+z\alpha^2& \alpha +z\alpha^2& \alpha +z\alpha & \alpha^2& 1+z\alpha
  \end{bmatrix}.
\end{equation}
Notice also that by Theorem~\ref{T-annihMat} the polynomial~$h'\in A[z;\widehat{\sigma}]$
is a generator polynomial of the dual code~$\cC^{\perp}$. A control
polynomial of that code is easily computed as
\[
  g':=\wh{g}{\sigma}=1+x^2+\alpha^2x^3+\alpha^2x^4+z(1+x^2+\alpha^2x^3+\alpha^2x^4)
\]
Altogether we have
\[
   \p(\cC)=\lideal{g}=\lcirc{\rideal{h}}\subseteq\Azs, \quad
   \p(\cC^{\perp})=\lideal{h'}=\lcirc{\rideal{g'}}\subseteq A[z;\widehat{\sigma}].
\]
It is worth mentioning that the code~$\cC$ has free distance equal to~$8$.
This is optimal among all codes
with the same parameters (length~$n=5$, dimension~$k=2$, complexity~$\delta=2$,
memory~$m=1$, and field size~$|\F|=q=4$) according to the Heller bound
(see \cite[p.~132]{JoZi99} for the binary case)
\[
  \dfree\leq\min\left\{\left\lfloor
     \frac{n(m+i)q^{k(m+i)-\delta-1}(q-1)}{q^{k(m+i)-\delta}-1}\right\rfloor\,\Bigg|\;
            i\in\N\right\}
\]
(the memory is the largest row degree appearing in a minimal generator
matrix in the sense of Definition~\ref{D-minimal} below).
The free distance of the dual is~$5$, attained by the
constant codeword
$v:=(\alpha,\alpha,\alpha,\alpha,\alpha)=(1,\alpha,1,0,0)\trans{\cMsigma{h}}$.
The dual code can also be regarded as optimal among all codes with the same parameters,
since each code with complexity~$2$
and dimension~$3$ has to contain a constant codeword.
This follows from the existence
of minimal generator matrices in the sense of Definition~\ref{D-minimal}
and the alternative characterizations given in~\cite[p.~495]{Fo75}.
\mbox{}\hfill$\Box$
\end{exa}

Next we will investigate the dimension of a $\sigma$-CCC, i.~e. the rank
of a $\sigma$-circulant, in terms of a given generator polynomial.
For a classical circulant the rank of~$M_g$ can be read off from the
polynomial~$g\in A$ via the formula given in Lemma~\ref{R-MgMh}(d).
Furthermore, the result shows how to cut out a rectangular generator
matrix of full rank from the square singular circulant~$M_g$.
For \scirc s these results are not true in this generality.
For instance, the matrix $\cMsigma{h}$ above is
basic of rank~$3$, hence $\im\cMsigma{h}$ is a $3$-dimensional
$\sigma$-CCC, but the first~$3$ rows do not form a generator matrix
of that code, since one can show that
$z^2+z\alpha+\alpha$ is a common factor of the full size minors of that
$3\times5$-matrix.

However, as we will show next, choosing a reduced generator~$g$ of the
ideal~$\p(\cC)$ always leads to a rectangular full rank generator matrix
of $\cC=\im\cMsigma{g}$ formed by the appropriate number of rows of the
\scirc. In order to prove this result we have to combine the techniques of
this section with the results and methods of the two foregoing sections.
It is quite advantageous to give some technicalities beforehand.
We will make use of the framework as in\eqnref{E-xn-1} --\eqnref{E-Kk}.
In particular, let $x^n-1=\pi_1\cdot\ldots\cdot\pi_r$ be the decomposition of
$x^n-1$ into its prime factors and for each $k=1,\ldots,r$, let $\ve{k}$ be
the irreducible idempotent associated with~$\pi_k$, i.~e.
$\ve{k}A=:K^{(k)}\cong\F[x]/{\ideal{\pi_k}}$.
Then $\pi_k\ve{k}=0$ when considered in~$A$ and
\begin{equation}\label{e-Pia1}
   \pi_k \mid a\ \text{ in } \F[x]\ \Longleftrightarrow \ a\ve{k}=0
\end{equation}
for all $a\in A$. For each $g\in\cR$ we define
\begin{equation}\label{e-Pia2}
     \pi_{(g)}:=\prod_{k\in T_g}\pi_k\in\F[x],
\end{equation}
where, as before, $T_g$ denotes the support of~$g$, see Notation~\ref{Nota1}.
Then $T_{\pi_{(g)}}\cap T_g=\varnothing$ and $T_{\pi_{(g)}}\cup T_g=\{1,\ldots,r\}$
by Equation\eqnref{e-Pia1} and thus
\begin{equation}\label{e-Pia3}
    \pi_{(g)}g=0\text{ in }\cR.
\end{equation}

In the case where $T_g=T_{g_0}$ one can alternatively express $\pi_{(g)}$
as $\pi_{(g)}=\frac{x^n-1}{\gcd(g_0,x^n-1)}$.

Now we are in a position to show

\begin{theo}\label{P-rank}
Let $g\in\cR$ be a reduced and nonzero polynomial and let
$\kappa:=\deg_x\pi_{(g)}$.
Then the family
\begin{equation}\label{E-family}
 g,xg,\dots ,x^{\kappa-1}g
\end{equation}
is a left $\F[z]$-basis of $\lideal{g}$.
Equivalently, $\rank\cMsigma{g}=\kappa$ and
the first $\kappa$ rows of $\cMsigma{g}$ form a full rank
generator matrix $G\in\F[z]^{\kappa\times n}$
of the $\sigma$-cyclic submodule $\cC:=\v\big(\lideal{g}\big)\subseteq\F[z]^n$.
Furthermore, if~$\cC$ is a code, i.~e. the matrix~$\cMsigma{g}$ is basic, then~$G$
is basic, too.
\end{theo}

\begin{proof}
Only the first part needs proof.
In order to establish left $\F[z]$-independence of the family in\eqnref{E-family}
let $u_i=\sum_{\nu=0}^{d_i}z^{\nu}u_{i\nu}$, where $u_{i\nu}\in\F$ and $i=0,\ldots,\kappa-1$,
and suppose $\sum_{i=0}^{\kappa-1}u_ix^ig=0$.
Accepting possible zero-coefficients, we may assume
$d_i=d$ for $0\le i\le \kappa-1 $. Letting $f_{\nu}:=\sum_{i=0}^{\kappa-1}u_{i\nu}x^i \in A$
and $f:=\sum_{\nu=0}^dz^{\nu} f_{\nu}$ we obtain
$$0=fg=fg^{(1)}+\cdots +fg^{(r)}$$
and by Lemma~\ref{rules}(c) and~\ref{norming}(d) we conclude that
$f\ve{k}=0$ for $k\in T_g$.
The definition of~$f$ shows that then also $f_{\nu}\ve{k}=0$ for
all $0\le \nu \le d$ and all $k\in T_g$.
Using\eqnref{e-Pia1} we obtain
$\pi_k \mid f_{\nu} \text{ in } \F[x]$ for all $k\in T_g$ and thus
$\pi_{(g)}  \mid  f_{\nu} \ \text{ in } \F[x]$.
Since $\deg_x (\pi_{(g)})=\kappa>\deg_x f_{\nu}$, the latter implies $f_{\nu}=0$
for all $\nu=0,\ldots,d$.
But then also $u_0=\ldots=u_{\kappa-1}=0$, showing the independence of the
given family.
\\
It remains to show that for $\kappa{'}\ge \kappa$ the polynomial $x^{\kappa{'}}g$
can be generated with coefficients from~$\F[z]$ by the family\eqnref{E-family}.
This is true even with coefficients from~$\F$ as can be deduced
recursively by using\eqnref{e-Pia3}.
\end{proof}

One should observe that a constant polynomial, i.~e.~$g\in A$, is
always reduced and in this case
$\pi_{(g)}=\frac{x^n-1}{\gcd(x^n-1,g)}$, see\eqnref{e-Pia1}.
Hence Theorem~\ref{P-rank} provides a generalization of the rank
formula for classical circulants given in Lemma~\ref{R-MgMh}(d).

The last part of the proof above shows that even for non-reduced
polynomials~$g$ the family in\eqnref{E-family} is a generating
system of the left $\F[z]$-module $\lideal{g}$.
However, in this case the family need not be independent, or
equivalently, $\kappa$ might be strictly bigger than $\rank\cMsigma{g}$.
We will show an example below in part~(3).

\begin{exa}\label{E-rank}
Let us reconsider Example~\ref{E-sigma.tau2} along with the various
representations.
\begin{arabiclist}
\item The polynomial~$g$ is reduced since
      $\varrho(g)=\ve{3}(\alpha x+1)+z\ve{2}(\alpha^2x+\alpha^2)=\ve{3}\varrho(g)$.
      Normalization of this polynomial has been performed in
      Example~\ref{E-norming}.
      As stated in Example~\ref{E-sigma.tau2}, the associated \scirc\
      has rank~$2$ which is also in accordance with the theorem
      above since $\pi_{(g)}=\pi_3=x^2+\alpha^2x+1$.
      Furthermore, as stated in the theorem, the first two rows of
      $\cMsigma{g}$ form a generator matrix of the code
      $\cC=\v\big(\lideal{g}\big)$, which can also be checked directly.
\item The dual code is given by
      $\cC^{\perp}=\im\cM^{\widehat{\sigma}}(h')$ where~$h'$ is as
      in\eqnref{e-hprime}. Since~$h'$ was the output of the reduction
      algorithm, it is reduced and thus Theorem~\ref{P-rank} is
      applicable again. As can be seen from\eqnref{e-h2prime} we now have
      $\pi_{(h)}=\pi_1\pi_2$, thus $\kappa=3$ telling us that the first
      three rows of~$\cM^{\widehat{\sigma}}(h')$ form a (basic) matrix of
      rank~$3$.
\item Let us also consider the code $\cC':=\im\cMsigma{h}\subseteq\F[z]^n$,
      where~$h$ is as in\eqnref{e-hpoly}. In this case the polynomial~$h$ is not
      reduced as one can see from $\varrho(h)=\ve{2}+\ve{1}+z\ve{2}$.
      The matrix $\cMsigma{h}$ is basic of rank~$3$ (see
      Lemma~\ref{L-annih}) but the first three rows do not span the
      code~$\cC'$. Reduction of~$h$ leads to the polynomial~$\tilde{h}$ where
      $\varrho(\tilde{h})=\ve{2}+\ve{1}$.
      Since~$\tilde{h}\in A$, we now get that $\cMsigma{\tilde{h}}=M_{\tilde{h}}$
      is a classical circulant and the code~$\cC'$ a $3$-dimensional
      cyclic block code.
      Let us compare this with Theorem~\ref{P-rank}.
      Despite the non-reducedness of the polynomial we can calculate the
      polynomial~$\pi_{(h)}$ and obtain $\pi_{(h)}=\pi_1\pi_2\pi_3=x^5-1$.
      Thus $\kappa=5>\rank\cMsigma{h}$ and the family in\eqnref{E-family}
      is not $\F[z]$-linearly independent, but certainly an
      $\F[z]$-generating set of~$\lideal{h}$.
      On the other hand, the reduced polynomial~$\tilde{h}$ satisfies
      $\pi_{(\tilde{h})}=\pi_1\pi_2$, thus $\kappa=3$ in accordance
      with $\rank\cMsigma{h}=\rank\cMsigma{\tilde{h}}=3$.
      \mbox{}\hfill$\Box$
\end{arabiclist}
\end{exa}


At this point it might also be interesting to know whether one can tell
from a given polynomial~$g\in\cR$ if $\cMsigma{g}$ is basic, in other
words, if $\v(\lideal{g})$ defines a code.
In case of a reduced polynomial~$g$ this can be characterized as follows.

\begin{prop}\label{P-polybasic}
Let $g\in\cR$ be a nonzero reduced polynomial with $z$-free term~$g_0$.
Then
\begin{align}
  \cMsigma{g}\text{ is basic }&\Longleftrightarrow
     \lideal{\wh{g}{\sigma}}=\lideal{\widehat{g_0}}
         \text{ in }A[z;\widehat{\sigma}] \nonumber \\
  &\Longleftrightarrow u\wh{g}{\sigma}=\widehat{g_0}\in A[z;\widehat{\sigma}]
      \text{ for some unit u in }A[z;\widehat{\sigma}]\nonumber \\
  &\Longleftrightarrow gv=g_0\text{ for some unit $v$ in }\Azs \label{e-gvv0}\\
  &\Longleftrightarrow \rideal{g}=\rideal{g_0}. \nonumber
\end{align}
In other words, $\cMsigma{g}$ is basic if and only if
$\trans{\!\cMsigma{g}}$ generates the cyclic block code
$\im\trans{M_{g_0}}$.
\end{prop}

One should note that the second equivalence says that $\wh{g}{\sigma}$ is left reducible
to the constant $\widehat{g_0}$.
It can be shown by examples, that the equivalence is not true if~$g$ is not
reduced.

\begin{proof}
Let $\pi=\pi_{(g)}$ be as in\eqnref{e-Pia2}. Then
$\cMsigma{\pi}\cMsigma{g}=M_{\pi}\cMsigma{g}=0$
and by Theorem~\ref{P-rank} (see also Lemma~\ref{R-MgMh}(d)) we have
$\rank\cMsigma{g}=n-\rank M_{\pi}$.
Therefore and upon using Proposition~\ref{P-directsummand} and
Lemma~\ref{L-annih} we obtain
\[
 \cMsigma{g}\text{ basic }\Longleftrightarrow
 \im\trans{\!\cMsigma{g}}=\ker\trans{M_{\pi}}
 \Longleftrightarrow
 \im\cM^{\widehat{\sigma}}(\wh{g}{\sigma})=\ker M_{\widehat{\pi}}
 \Longleftrightarrow
 \lideal{\wh{g}{\sigma}}=\lcirc{\rideal{\widehat{\pi}}}.
\]
Since $\widehat{\pi}\in A$ we have
$\lcirc{\rideal{\widehat{\pi}}}=\lideal{\widehat{a}}$ for
$\widehat{a}=\frac{x^n-1}{\gcd(x^n-1,\widehat{\pi})}\in A$.
Hence $\cMsigma{g}$ is basic if and only if
$\wh{g}{\sigma}$ can be left reduced to the constant~$\widehat{a}$.
By Corollary~\ref{C-redu} this is equivalent to
the existence of some unit $u\in A[z;\widehat{\sigma}]$ such that
$u\wh{g}{\sigma}=\widehat{a}$.
Since the $z$-free term~$u_0$ of~$u$ is a unit in~$A$ and
the $z$-free term of $\wh{g}{\sigma}$ is given by $\widehat{g_0}$, we obtain
the identity $u_0\widehat{g_0}=\widehat{a}$ and without restriction we may assume
$\widehat{a}=\widehat{g_0}$.
This yields the desired result.
\end{proof}

The proposition above has an interesting consequence.

\begin{cor}\label{C-commut}
Let $g,\,h\in\cR$ such that $\rcirc{\lideal{g}}=\rideal{h}$ and
$g\in\Azs$ and $\wh{h}{\sigma}\in A[z;\widehat{\sigma}]$ are both left reduced, which
can be assumed without restriction.
Furthermore, assume that~$g$ generates a code, thus $\cMsigma{g}$ is basic.
Then $hg=0$ and even
\[
    \rcirc{\lideal{h}}=\rideal{g}.
\]
In particular, the identity $\cMsigma{g}\cMsigma{h}=0$ implies
that also $\cMsigma{h}\cMsigma{g}=0$.
\end{cor}

Again, the result is not true if any of the polynomials is not reduced.

\begin{proof}
First notice that $\cMsigma{h}$ is basic by assumption, see Lemma~\ref{L-annih}.
Thus, we may apply Proposition~\ref{P-polybasic} to the polynomials~$g$
and~$\wh{h}{\sigma}$ in their respective Piret algebras and obtain
$gu=g_0$ and $\wh{h}{\sigma}v=\widehat{h_0}$ for some units $u\in\Azs$
and $v\in A[z;\widehat{\sigma}]$.
Then $gh=0$ implies $0=g_0h_0=h_0g_0$, since~$A$ is commutative, and thus
$0=\wh{v}{\widehat{\sigma}}hgu$, after applying the anti-isomorphism
$\wh{\ }{\widehat{\sigma}}$.
Cancellation of the units yields $hg=0$ and thus
$\im\cMsigma{h}\subseteq\ker\cMsigma{g}$.
Furthermore, from Lemma~\ref{L-annih} we know that $\rank\cMsigma{g}=n-\rank\cMsigma{h}$
and since $\cMsigma{h}$ is basic we may apply Proposition~\ref{P-directsummand}(7)
in order to get $\im\cMsigma{h}=\ker\cMsigma{g}$.
Then Lemma~\ref{L-annih}(2) completes the proof.
\end{proof}

As a by-product, Proposition~\ref{P-polybasic} gives us an alternative proof of
Proposition~\ref{P-trivialCCC} since in the case where $\sigma=\text{id}$, the ring
$\Azs$ is commutative
and therefore\eqnref{e-gvv0} is the same as $vg=g_0$ so that,
consequently, the corresponding left ideal has a constant generator.

Finally, it remains to discuss the important issue of minimal
generator matrices. In convolutional coding theory one is mainly
interested in minimal encoding matrices since they have, by
definition, minimum possible row degrees, so that, as a
consequence, their canonical linear shift realization need the
minimum number of memory elements; for details
see~\cite[Sec.~2.7]{JoZi99}. The row Hermite form of a polynomial
matrix usually tends to have artificially high degrees in its
entries and therefore is not minimal. The following definition is
adapted to our purposes. More common but equivalent definitions
can also be found e.~g. in~\cite[p.~495]{Fo75}.

\begin{defi}\label{D-minimal}
Let $M\in\F[z]^{m\times n}$ be a matrix with rows $w_1,\dots ,w_m\in\F[z]^n$ and
$\rank_{\F[z]}M=m$.
The leading $z$-coefficient vector of~$w_i$ will be denoted by $\lcz{w_i}\in\F^n$.
The matrix~$M$ is called (row-) minimal if its (row-) leading coefficient matrix
\[
    L(M):=
    \begin{bmatrix}
         \lcz{w_1}\\ \lcz{w_2}\\ \vdots \\ \lcz{w_m}
    \end{bmatrix}\in \F^{m\times n}
\]
satisfies $\rank_{\F}L(M)=m$.
\end{defi}

It can easily be seen via some examples that the full rank generator matrix of a $\sigma$-CCC
as constructed in Theorem~\ref{P-rank} in general is not minimal.
This is, for instance, the case for the matrix~$\widehat{G}$  formed by the
first three rows of~$\cM^{\widehat{\sigma}}(h')=\trans{\cMsigma{h}}$ in
Example~\ref{E-sigma.tau2}.
The matrix $\widehat{G}$ is a basic generator matrix of the dual code~$\cC^{\perp}$, but
not minimal.

We will now show, how one can obtain a minimal generator matrix
by extracting the appropriate number of first rows of the circulants associated
with the components of a reduced generator polynomial.

\begin{theo}\label{T-minimal}
\begin{alphalist}
\item Let $g^{(k)}\in\ve{k}\cR$ be non-zero and let $\pi_k$ be the prime
      divisor of $x^n-1$ corresponding to $\ve{k}$. Put $\kappa_k:=\deg_x{\pi_k}$.
      Then the matrix
      \[
         G_k:=\begin{bmatrix}
               \v(g^{(k)})\\ \v (xg^{(k)})\\ \vdots \\ \v (x^{\kappa_k-1}g^{(k)})
              \end{bmatrix}\in\F[z]^{\kappa_k\times n}
      \]
      formed by the first~$\kappa_k$ rows of $\cMsigma{g^{(k)}}$
      is a minimal generator matrix for the $\F[z]$-module
      $\v\big(\lideal{g^{(k)}}\big)\subseteq \F[z]^n$.
\item Let $g\in \cR$ be non-zero and left reduced.
      Suppose $T_g= \{k_1,\dots ,k_t\}$ and put
      $\kappa_{k_{\nu}}:=\deg_x{\pi_{k_{\nu}}}$ for $1\le \nu \le t$
      and $\kappa:=\sum_{\nu=1}^t\kappa_{k_{\nu}}$,
      \[
         G:=\begin{bmatrix}
           G_{k_1}\\ \vdots\\G_{k_t}
            \end{bmatrix}\in\F[z]^{\kappa\times n}
         \text{\quad and }\quad
          G_{k_{\nu}}:=
          \left[
          \begin{array}{c}
          \v (g^{(k_{\nu})})\\ \vdots\\ \v (x^{\kappa_{k_{\nu}-1}}g^{(k_{\nu})})
          \end{array}
          \right].
      \]
      Then $G$ is a minimal generator matrix for
      the $\F[z]$-module $\v\big(\lideal{g}\big)\subseteq \F[z]^n$.
\end{alphalist}
\end{theo}

\begin{proof}
(a)
Let $\deg_z{g^{(k)}}=d_k$ and denote the leading $z$-coefficient of~$g^{(k)}$
by~$g_{d_k}$, which then is nonzero.
From Theorem~\ref{P-rank} we know that
$g^{(k)},\dots ,x^{\kappa_k-1}g^{(k)}$ is a left $\F[z]$-basis
for $\lideal{g^{(k)}}$ and that~$G_k$
is a full rank generator matrix of $\v\big(\lideal{g^{(k)}}\big)$.
It remains to check minimality.
Note that for all $i=0,\ldots,\kappa_k-1$ the leading $z$-coefficient of the
polynomial $x^ig^{(k)}$ is given by~$\sigma^{d_k}(x^i)g_{d_k}$.
Therefore, the leading coefficient matrix of~$G_k$ is
\[
  L(G_k)=
   \begin{bmatrix}
      \v(g^{(k)})\\ \v(\sigma^{d_k}(x)g_{d_k})\\ \vdots\\
      \v(\sigma^{d_k}(x^{\kappa_k-1})g_{d_k})
   \end{bmatrix}
\]
and we have to show that its rank is equal to~$\kappa_k$.
To this end suppose
\[
  \sum_{i=0}^{\kappa_k-1}c_i\v (\sigma^{d_k} (x^i)g_{d_k})=0 \text{ for some }
  c_0,\dots , c_{\kappa_k-1}\in\F.
\]
Then we compute
\begin{align*}
   0=\v \big(\sum_{i=0}^{\kappa_k-1}c_i \sigma^{d_k} (x^i)g_{d_k}\big)
   =\v\bigg(\sigma^{d_k}
      \Big(\Big(\sum_{i=0}^{\kappa_k-1}c_ix^i\Big)\sigma^{-d_k}(g_{d_k})\Big)\bigg).
\end{align*}
Since $\sigma^{-d_k} (g_{d_k})$ is from $\ve{k}A$ and,
of course, also nonzero and since $\ve{k}$ is idempotent, we may
use\eqnref{E-idemzero} and conclude
\[
0=\sum_{i=0}^{\kappa_k-1}c_i\big(x\ve{k}\big)^i.
\]
Since $c_i\in\F$, this equation takes place in the field $K^{(k)}=\ve{k}A$, and
$\kappa_k\leq n$ implies $c_0=\dots=c_{\kappa_k-1}=0$.
\\
(b)
By (a) we know that the family
$g^{(k_{\nu})},\dots ,x^{\kappa_{k_{\nu}}-1}g^{(k_{\nu})}$ generates
$\lideal{g^{(k_{\nu})}}$ for $1\le \nu\le t$. Therefore the~$t$ families
together generate the $\F[z]$-left module $\lideal{g^{(k_1)},\ldots,g^{(k_t)}}=\lideal{g}$.
Reducedness of~$g$ and Theorem~\ref{P-rank} imply that the $\F[z]$-rank of $\lideal{g}$
is~$\kappa$.
Recalling that $xg^{(k_{\nu})}=(xg)^{(k_{\nu})}$ we therefore see, that
the entire family  $\big(x^ig^{(k_{\nu})}\big)_{0\le i\le\kappa_{k_{\nu}} -1,\,1\le \nu\le t}$
is $\F[z]$-linearly independent.
This guarantees that~$G$ has full rank and it remains to consider the
leading coefficient matrix~$L(G)$. This time we have
\[
   L(G)=\begin{bmatrix}
        \v (g_{d_1})\\ \vdots\\
        \v (\sigma^{d_1}(x^{\kappa_{k_1}-1}) g_{d_1})\\
        \vdots\\ \vdots\\
        \v (g_{d_t})\\ \vdots\\
        \v (\sigma^{d_t}(x^{\kappa_{k_t}-1}) g_{d_t})
      \end{bmatrix},
\]
where $d_{\nu}=\deg_z g^{(k_{\nu})}$ and
$g_{d_{\nu}}\ne0$ is the leading $z$-coefficient of~$g^{(k_{\nu})}$.
Suppose now $c\,L(G)=0$ for some vector
$c=(c_{10},\dots ,c_{1\,\kappa_{k_1-1}},\dots ,
   c_{t0},\dots ,c_{t\,\kappa_{k_t-1}})\in\F^{\kappa}$.
Then we conclude as in~(a)
\begin{equation}\label{E-minim}
    0=\sum_{\nu =1}^t\sum_{i=0}^{\kappa_{k_{\nu}}-1}
    c_{\nu i}\sigma^{d_{\nu}}(x^i)g_{d_{\nu}}=
    \sum_{\nu=1}^t\bigg(\sigma^{d_{\nu}}
    \bigg(\sum_{i=0}^{\kappa_{k_{\nu}}-1}c_{\nu i}\big(x\ve{k_{\nu}}\big)^i\bigg)
    g_{d_{\nu}}\underbrace{\sigma^{d_{\nu}}(\ve{k_{\nu}})}_{\text{idempotent}}\bigg).
\end{equation}
Since $g$ is reduced, no two of the idempotents
$\sigma^{d_1}(\ve{k_1}), \dots , \sigma^{d_t}(\ve{k_t})$
can be equal. Therefore Equation\eqnref{E-minim} implies
\[
    \sum_{i=0}^{\kappa_{k_{\nu}}-1}c_{\nu i}(x\ve{k_{\nu}})^i=0 \
    \text{ for all }\ 1\le\nu\le t.
\]
Just like in~(a) we conclude $c=0$ and the matrix $L(G)$ has full rank.
\end{proof}

\begin{exa}\label{E-minbasis}
Consider again Example~\ref{E-sigma.tau2}.
In Example~\ref{E-rank} we saw already that $\varrho(g)=\ve{3}\varrho(g)$ is
reduced. According to the theorem above the first two rows of $\cMsigma{g}$ form a minimal
basic generator matrix of the code $\cC=\v(\lideal{g})$, which can also
be seen directly from the matrix given in Example~\ref{E-sigma.tau}.
Furthermore, the first three rows of
$\trans{\cMsigma{h}}=\cM^{\widehat{\sigma}}(h')$ form a basic generator
matrix of the code $\cC^{\perp}$.
But as is easily seen, the matrix is not minimal.
According to the theorem above and the representation\eqnref{e-h2prime}
we have to combine the first row of
$\cM^{\widehat{\sigma}}\big(\varrho^{-1}(\ve{1}h')\big)$ and the first two
rows of
$\cM^{\widehat{\sigma}}\big(\varrho^{-1}(\ve{2}h')\big)$
in order to get a minimal basic generator matrix of the code $\cC^{\perp}$.
This leads to the matrix
$$
   \begin{bmatrix} 1&1&1&1&1\\0&
     \alpha^2 z+\alpha&\alpha z+\alpha^2&\alpha z+\alpha^2&\alpha+\alpha^2z\\
     \alpha z+\alpha&\alpha z&\alpha^2z+\alpha&\alpha^2&\alpha^2z+\alpha^2
   \end{bmatrix}.
   \eqno\Box
$$
\end{exa}

The last theorem allows for a formula for the complexity of a
$\sigma$-cyclic code in terms of a reduced generator. The key
point is that the complexity, as defined in
Definition~\ref{D-CC}(4), can be computed much easier if a minimal
generator matrix is available. Indeed, it is known
from~\cite[p.~495]{Fo75}, see also~\cite[Sect.~3]{McE98}, that if
$\cC=\im G$ where $G\in\F[z]^{k\times n}$ is a minimal matrix with
rows $G_1,\ldots,G_k$, then the complexity is given by
$\delta=\sum_{i=1}^k \deg_z G_i$. Using
Theorem~\ref{T-minimal}(b), this immediately implies

\begin{cor}\label{C-compl}
Let $g\in\cR$ be a reduced polynomial such that $\cMsigma{g}$ is basic and let
$\cC:=\im\cMsigma{g}$ be the $\sigma$-cyclic code generated by~$g$.
Then the complexity of~$\cC$ is given by
\[
   \delta=\sum_{i\in T_g} \deg_x\pi_i\deg_z g^{(i)},
\]
where, again, $\pi$ is the prime divisor of $x^n-1$ corresponding to~$\ve{i}$.
\end{cor}

It remains to present the

{\sc Proof of Proposition~\ref{P-blockcode}:}
(a)
``$\Rightarrow$'': Let $\sigma (\ve{k})=\ve{l}\ne\ve{k}$
and put $g:=z\ve{l}+\ve{k}$. Note that $g=\ve{k}g$. We
claim that $g$ generates a left ideal
$\cJ=\lideal{g}$ corresponding to a $\sigma$-CCC, which cannot be
generated by a constant matrix. In order to prove this it suffices
to show that, firstly, $\cJ$ has no constant generator and
that, secondly,~$\cJ$ is a direct summand as a left $\F[z]$-submodule
of~$\cR$, see Observation~\ref{O-attrib}.
A constant generator necessarily would be $\ve{k}$ up to a unit from~$A$.
Hence assume $\ve{k}=vg$ for some $v\in\cR$.
Comparing like powers of~$z$ in the equation $\ve{k}=vg=v\ve{k}(z+1)$
shows that this is not possible since the leading coefficient of $z+1$ is
a unit in~$A$.
Therefore $\cJ$ has no constant generator.
As for the direct summand property, assume $fu=vg\in\lideal{g}=\cJ$ for some $f\in\F[z]$
and $u,\,v\in\cR$. But then also $vg\ve{k}=v\ve{k}=fu\ve{k}$ and thus
$fu=vg=v\ve{k}g=fu\ve{k}g$. But the latter implies $u\in\cJ$,
since $f\in\F[z]$, not being a zero divisor in $\cR=\Azs$, can be cancelled.
Hence~$\cJ$ is a direct summand of $\cR$.
\\
``$\Leftarrow$'':
The assumption $\sigma (K^{(k)})=K^{(k)}$ for $1\le k\le r$ can be
rephrased as $\sigma (\ve{k})=\ve{k}$ for  $1\le k\le r$.
This in turn implies that all idempotents are lying in the
center of~$\cR$ i.e. $\ve{k}g=g\ve{k}$ for all $1\le k\le r$ and all
$g\in\cR$.
Now let $\cC$ be a $\sigma$-CCC and
$\cJ=\p (\cC)$ be the corresponding left ideal. We have to
show that $\cJ$ has a constant generator polynomial.
Since $\cC$ is delay-free we know from Corollary~\ref{C-redu}
that $\cJ=\lideal{g}$ for some reduced polynomial $g\in\cR$
which also satisfies\eqnref{E-existence}.
Define now $\varepsilon:=\sum_{k\in T_g} \ve{k}$. Then
$\varepsilon g=g=g \varepsilon$
and as a consequence $\cJ=\lideal{g}\subseteq \lideal{\varepsilon}$.
The polynomials~$g$ and $\varepsilon$ are both reduced and
satisfy $T_g=T_{\varepsilon}$. Therefore
Theorem~\ref{P-rank} yields that~$\cJ$ and $\lideal{\varepsilon}$ have the same rank as
$\F[z]$-submodules of $\cR$. Since~$\cJ=\p(\cC)$ is a direct summand
it follows $\cJ=\lideal{\varepsilon}$, see Proposition~\ref{P-directsummand}(7),
showing that~$\cJ$ has a constant generator.
\\
(b) can be shown with exactly the same line of arguments as
in~``$\Leftarrow$''of~(a).
\mbox{}\hfill$\Box$

\section{Future research topics}\label{S-outlook} \setcounter{equation}{0}

In this paper we made an effort to broaden the mathematical basis
for a thorough investigation of $\sigma$-cyclic convolutional
codes and their potential for coding.
Yet, many important questions of coding theory still have to be answered.
We hope that our contribution might serve as a basis for further investigations in this
direction and close the paper with a brief list of issues to be addressed in the future.
(a) In the paper~\cite{GSS02} we presented an infinite series of $1$-dimensional codes of
length~$2$ over $\F_3$ with increasing complexity. We also showed that the first codes
in this series have a pretty good distance. It would be worth knowing whether the
free distance of these codes tends to infinity for increasing complexity.
More generally, one might ask whether it is possible to construct families of
$\sigma$-cyclic codes with constant dimension and length over a fixed field and
with arbitrary large distance.
(b) Any convolutional code allows for other representations besides those via
generator and control matrices, see for instance~\cite[p.~1071]{McE98}
or~\cite{RSY96},
where a shift realization is translated into a description of the code as a first order
discrete-time dynamical system over the field~$\F$.
Is it possible to recover cyclic structure in this description? If so, can such a
description be used for the construction of good cyclic codes?
(c) One of the strengths of cyclic block codes is the relation between the zeros of
the generator polynomial and the distance of the code, leading to the design of powerful
codes like BCH-codes.
The central issue of the theory of CCC's is certainly the investigation of the distance of these
codes in terms of a generator or control polynomial or other data determining the code.
Any algebraic result in this direction would improve the theory of CCC's.
(d) The other main advantage of cyclic block codes is their potential for decoding.
Does the additional structure of CCC's, beyond the $\F[z]$-module structure, also allow for an
algebraic decoding algorithm, that is, an algorithm where decoding is not obtained via a search
algorithm but rather via an algebraic computation based on the received word?
A positive answer would certainly be a breakthrough in the theory of convolutional codes.

%
\bibliographystyle{abbrv}
\bibliography{literatureAK,literatureLZ}
\end{document}